\documentclass[12pt,a4paper]{article}

\usepackage{amsmath}    
\usepackage{amssymb}
\usepackage{graphicx}   
\usepackage{verbatim}   
\usepackage{color}      

\usepackage{hyperref}   

\usepackage{enumerate}
\usepackage{mathrsfs}
\usepackage{empheq}
\usepackage{bbold}
\usepackage{empheq}
\usepackage{xspace}
\usepackage{pdfpages}

\usepackage[final]{showlabels}

\usepackage{amsthm}
\newtheorem{definition}{Definition}[section]
\newtheorem{lemma}{Lemma}[section]

\newtheorem{theorem}{Theorem}[section]
\newtheorem{example}{Example}[section]

\newcommand{\del}{\partial}
\renewcommand{\theta}{\vartheta}
\renewcommand{\phi}{\varphi}

\newcommand{\dd}{\mathrm{d}}

\newcommand{\const}{\mathrm{const}}

\newcommand{\admis}{\ensuremath{\mathcal A}\xspace}

\renewcommand{\title}{A well-balanced Active Flux method for the shallow water equations with wetting and drying}

\newcommand{\authorOne}{Wasilij Barsukow\footnote{Bordeaux Institute of Mathematics, Bordeaux University and CNRS/UMR5251, Talence, 33405 France}}
\newcommand{\authorTwo}{Jonas\,P. Berberich\footnote{Würzburg University, Emil-Fischer-Strasse 40, Würzburg, 97074 Germany}}

\begin{document}

\begin{center} \Large
\title

\vspace{1cm}

\date{}
\normalsize

\authorOne, \authorTwo
\end{center}

\begin{abstract}

Active Flux is a third order accurate numerical method which evolves cell averages and point values at cell interfaces independently. It naturally uses a continuous reconstruction, but is stable when applied to hyperbolic problems. In this work, the Active Flux method is extended for the first time to a nonlinear hyperbolic system of balance laws, namely to the shallow water equations with bottom topography. We demonstrate how to achieve an Active Flux method that is well-balanced, positivity preserving, and allows for dry states in one spatial dimension. Because of the continuous reconstruction all these properties are achieved using new approaches. To maintain third order accuracy, we also propose a novel high-order approximate evolution operator for the update of the point values. A variety of test problems demonstrates the good performance of the method even in presence of shocks.



\end{abstract}

\section{Introduction}

The shallow water equations, derived by Adhémar Barré de Saint-Venant (\cite{saintvenant71}),
\begin{align}
 \del_t h + \del_x m &= 0 & h &: \mathbb R^+_0 \times \mathbb R \to \mathbb R^+ \label{eq:shallowwaterh}\\
 \del_t m + \del_x\left( \frac{m^2}{h} + \frac12 g h^2 \right) &= - g h \del_x b & m&: \mathbb R^+_0 \times \mathbb R \to \mathbb R, \,\, g \in \mathbb R\label{eq:shallowwaterm}
\end{align}
model the inviscid dynamics of a free water surface under the influence of gravity. Here, $h$ is the water height, $v = m/h$ the speed of the water (height-averaged), $b: \mathbb R \to \mathbb R $ the topography of the bottom and $g > 0$ the downwards gravitational acceleration. \eqref{eq:shallowwaterh}--\eqref{eq:shallowwaterm} is a hyperbolic system of balance laws. This system is useful to describe water of depth much smaller than the amplitude and wavelength of typical surface disturbances, hence the name. In relevant applications, any numerical method aiming at solving \eqref{eq:shallowwaterh}--\eqref{eq:shallowwaterm} needs to be able to deal with inundation, i.e.\ with the appearance and disappearance of dry states ($h=0$). It also needs to account for the possibility of stationary solutions, where the flux divergence balances the source term. Although the most general stationary solution may be a moving equilibrium ($m\neq 0$), here we focus on a particularly important special case, the so-called \emph{lake at rest}
\begin{align}
 h + b &= \const & m&= 0
\end{align}
In simulations, exact preservation of this stationary solution is required in order for the dynamics of small waves on top of such a lake not to be spoiled by spurious waves generated by an inaccurate numerical method. Numerical methods that preserve (discretizations of particular) stationary solutions exactly are called \emph{well-balanced}. 

The Active Flux method is an extension of the finite volume method and allows to stably solve hyperbolic problems with explicit time integration. It has a number of advantages over standard finite volume schemes. It incorporates additional (pointwise) degrees of freedom, which are evolved independently from the cell average. Active Flux thus is a \emph{high order method with a compact stencil}. These additional point values are distributed at the cell interfaces, which has two consequences. First, the numerical flux needed for the average update is immediately available, and second, they enforce a continuous reconstruction. Thus, Active Flux does not require Riemann solvers. This might be at least partially the reason why Active Flux is naturally structure preserving: it is, for example, vorticity and stationarity preserving for multi-dimensional acoustics (see \cite{barsukow18activeflux}) and also naturally well-balanced for acoustics with gravity (see \cite{barsukow19activefluxsource}). The time integration does not employ the method of lines and thus does not require to use e.g.\ Runge-Kutta methods. Active Flux thus is a promising candidate for combining high order of accuracy, upwinding (i.e.\ stability) and structure preservation.

To study and exploit its structure preserving capabilities, in this paper we for the first time propose an Active Flux method for the shallow water equations. Due to the inclusion of point values and a continuous reconstruction, it requires a different design in many aspects. In particular, in this paper, we present a novel positivity preserving reconstruction and an evolution procedure for Active Flux which allow for dry states. We also show how to obtain a well-balanced property for our method.

To remain in concordance with the continuous reconstruction of the conserved quantities, in this paper we choose to project the (given) bottom topography onto a globally continuous, piecewise parabolic function. This can be considered as a high order extension of the globally continuous, piecewise linear bottom topography considered in (\cite{ricchiuto09,bollermann13,kurganov18}). Many existing numerical methods use a bottom topography with jumps (for instance piecewise constant) (e.g.\ \cite{Bermudez1994,gallouet03,Audusse2004,castro06}), a natural choice when using Riemann solvers, but {requiring special effort in order to give sense to the derivative of the bottom topography $b$} (e.g.\ \cite{pares04}). {We believe that in practical applications, the continuity of the bottom topography is not a restriction.}

The above mentioned challenges of maintaining exactly the lake at rest and preventing negative water height have been addressed in a variety of ways in the literature (see e.g. \cite{xing14} for a review). \textbf{Well-balancing} strategies such as \cite{gosse01,Audusse2004} are designed to be used with Riemann solvers, and cannot be directly used for Active Flux. Other suggestions include kinetic schemes (e.g.\ \cite{xu02}) or formulations in equilibrium variables {(e.g.\ \cite{zhou01,kurganov02,rogers03,gallouet03,cheng19})}. We use this latter approach in cells far away from the shore, where we reconstruct $h+b$ instead of $h$. In principle, one would not even need to reconstruct in equilibrium  variables: With the natural reconstruction in Active Flux being piecewise parabolic and globally continuous, and given the piecewise parabolic, globally continuous $b$ one can reconstruct $h$ for a lake at rest exactly.  However, our limiting strategy employs other than parabolic reconstructions, at which point reconstruction of $h+b$ becomes inevitable. 

Equally, \textbf{positivity preserving\footnote{The word ``positivity'' is used throughout the paper to mean ``non-negativity'' in the context of dry states.}} limiters such as \cite{zhang10,xing11,castro19} only take into account the cell average and cannot be used if additionally the point values at cell boundaries are prescribed, as is the case for Active Flux. Here, we therefore for the first time extend the non-negative reconstruction from \cite{bollermann13}, which uses a piecewise linear, globally continuous bottom, to piecewise quadratic. {Although the high order of accuracy requires the consideration of more cases, the resulting reconstruction procedure remains computationally inexpensive, in particular because only a few cells are located close to the shore.} 

\begin{figure}
 \centering
 \includegraphics[width=0.9\textwidth]{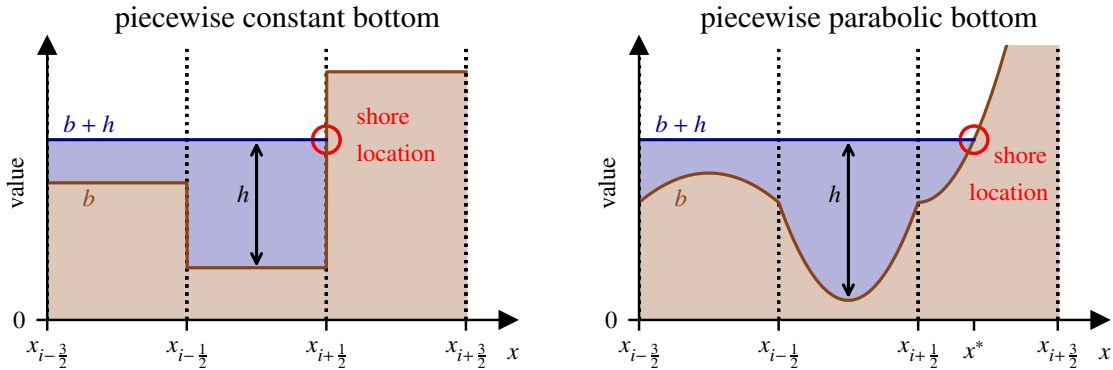}
 \caption{Different treatment of the shore for piecewise constant (\emph{left}, the shore is always located at a cell interface) and for globally continuous bottom (\emph{right}, the shore can be located anywhere inside the cell). In this paper, the latter approach is pursued. {Note that, in fact, for higher order versions of the left figure (e.g. for a piecewise linear bottom), the shore can be located anywhere inside the cell as well.}}
 \label{fig:constantbottomhaseasiershore}
\end{figure}

As described in \cite{kurganov18}, combining well-balancing and positivity preservation can be non-trivial if the bottom topography is globally continuous. Indeed, in this case there exist \textbf{half-wet cells} which contain the location of the shore. This is in contrast to a piecewise constant bottom, where cells are either dry or wet, and the location of the shore is always a cell interface (see Figure \ref{fig:constantbottomhaseasiershore}). Whereas in the lake-at-rest case, the water level is constant throughout a wet cell, it is only constant up to the shore in a cell that is partially dry. Thus, well-balancing of a lake at rest with a dry region now is more difficult, and is related to finding the location of the shore. {The same challenge, in fact, arises for e.g. piecewise linear reconstructions, and thus generally for high order methods. At the same time, the number of cells which contain a shore is very small such that their special treatment does not entail any real increase in computational time.} Here, we take inspiration from \cite{bollermann13,kurganov18} by choosing a reconstruction that is exact for the lake at rest and extending the approach to higher order. Such a reconstruction then makes sure that point values remain stationary for the lake at rest. 

Recall that Active Flux uses both point values and averages, and both need to be updated in a well-balanced fashion while avoiding negative water height. For the update of the averages, we observe that a particular source term quadrature ensures stationarity in the lake-at-rest case, which parallels findings made as early as in \cite{roe87,Bermudez1994,leveque98} (``reconstruction in piecewise equilibrium'', ``upwinding of the source term''). To ensure that the update of the cell averages is positivity preserving, we adopt the \emph{draining time-step technique} from \cite{bollermann11}.

\textbf{High-order} (i.e.\ higher than second order) well-balanced and positivity preserving methods have been suggested e.g. in \cite{xing05,xing06} based on WENO or Discontinuous Galerkin (DG) methods, and high-order methods which additionally allow for dry states are e.g. \cite{gallardo07,xing10,xing11,cozzolino12,castro19}. Active Flux is another way to achieve high order, and in this paper we apply it for the first time to the shallow water equations. It shares some similarity with DG methods, because it stores additional degrees of freedom inside the cell, thus keeping a compact stencil. However, these point values are placed at cell interfaces and shared between adjacent cells, which makes a huge difference to the approach of DG. The most significant difference to existing methods is the absence of discontinuities in both the reconstructed water height and the bottom topography. This greatly changes the entire approach, and, we believe, is often easier to handle. The maximum CFL number of Active Flux is optimal (i.e. $\mathrm{CFL}_\text{max} = 1$), while the shared degrees of freedom reduce the memory requirements. For an indication of favorable diffusion and dispersion characteristics of Active Flux in comparison to DG, see \cite{roe21}.

Because of the many aspects that need to be designed anew in order to fit the philosophy of Active Flux, this paper focuses on the one-dimensional case. This is largely due to the fact that multi-dimensional systems of hyperbolic PDEs do not generally have characteristics, but only characteristic cones, which makes the solution operators much more complicated -- even in the linear case (see e.g.\ \cite{barsukow17}). Currently, some approximate solution operators for multi-dimensional nonlinear systems have been suggested (e.g.\ \cite{fan17}), but many questions need to be answered before these operators can be adequately used in such an intricate setup as shallow water equations with dry states. Therefore, the multi-dimensional extension of the Active Flux method remains subject of future work.

The paper is organized as follows. Section \ref{sec:af} gives an overview of Active Flux and of its different ingredients. A globally continuous, non-negative reconstruction able to deal with dry and partially dry cells is presented in Section \ref{sec:bottomandrecon}. The update of the point values is described in Section \ref{sec:pointvalues}, in particular the way how well-balancing and the correct order of accuracy are achieved. Section \ref{sec:averages} discusses the update of cell averages, including the source term quadrature. Section \ref{sec:averages}, finally, shows a variety of numerical results that demonstrate the abilities of the new method.

\section{Overview of Active Flux} \label{sec:af}
\subsection{The algorithm}

The Active Flux method is an extension of the finite volume scheme, originating in \cite{vanleer77,eymann11,eymann13}. In the following, the general outline of Active Flux is given. Details specific to the shallow water equations are subject of later chapters. More thorough introductions to the approach of Active Flux are given in e.g. \cite{eymann13,kerkmann18,barsukow19activeflux}, to which the reader is referred.

\begin{figure}
	\centering
	\includegraphics[scale=1]{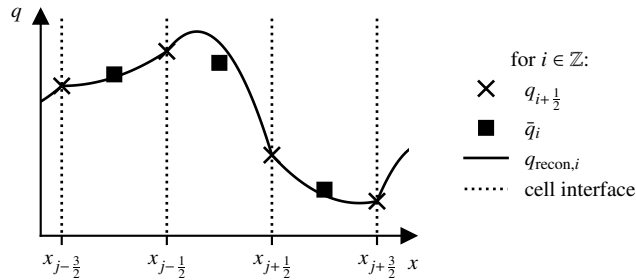}
	\caption{\label{fig:af_states}%
		The degrees of freedom in an Active Flux method are given by the interface values $q_{i+\frac12}$ (crosses) and the cell-average values $\bar q_i$ (squares). The continuous parabolic reconstruction $q_{{\rm recon},i}$ is uniquely determined by these degrees of freedom.
		}
\end{figure}

Active Flux incorporates additional degrees of freedom interpreted as point values. They are placed at cell interfaces and are shared by adjacent cells (see Figure \ref{fig:af_states}). In this paper only the one-dimensional situation is considered. The cell centers are denoted by $x_i$, the cell interfaces by $x_{i+\frac12}$, $i \in \mathbb Z$. $\Delta x$ shall denote the size of every cell, as the grid is assumed equidistant throughout the paper\footnote{However, due to the local nature of Active Flux, most of what is described immediately extends to non-equidistant grids, as long as the grid parameter varies continuously (i.e.\ $\Delta x_{i+1}/\Delta x_i \to 1$ upon refinement).}. Every cell therefore has access to its cell average $\bar q_i$ {(denoted by a bar and integer indices)} and two point values $q_{i\pm\frac12}$ {(denoted by indices in $\mathbb Z + \frac12$)} situated at the interfaces\footnote{We endow the average $\bar q_i^n$, the point values $q^n_{i+\frac12}$ and the reconstruction $q^n_{\text{recon},i}$ with temporal superscripts whenever necessary for the presentation, and -- in order to simplify notation -- refrain from making them explicit when no confusion is possible.}. Active Flux therefore is third order accurate.

A general $M \times M$ system of balance laws is given by
\begin{align}
 \del_t q + \del_x f(q) &= s(x, q) & q &: \mathbb R^+_0 \times \mathbb R \to \admis \label{eq:hypsystemgeneral}\\
  \nonumber &&f &: \admis \to \mathbb R^M\\
   \nonumber &&s &: \mathbb R \times \admis \to \mathbb R^M
\end{align}
where $\admis\subset \mathbb R^M$ denotes the set of admissible states for the system, e.g.\ $(h,b)\in\admis=\mathbb R^+\times\mathbb R$ for the shallow water system \eqref{eq:shallowwaterh}-\eqref{eq:shallowwaterm}.
The update of the cell average for \eqref{eq:hypsystemgeneral} is
\begin{align}
 \bar q_i^{n+1} = \bar q_i^n - \Delta t\frac{\hat f_{i+\frac12} - \hat f_{i-\frac12}}{\Delta x} + \Delta t \hat s_i \label{eq:generalafavgupdate}
\end{align}
where the discrete fluxes $\hat f_{i+\frac12}$ and source $\hat s_{i}$ are approximations of
\begin{align}
 \hat f_{i+\frac12} &\simeq \frac{1}{\Delta t} \int_0^{\Delta t} \dd t \, f(q(t, x_{i+\frac12})) \label{eq:fluxexact}\\
 \hat s_{i} &\simeq \frac{1}{\Delta t} \int_0^{\Delta t} \dd t \frac{1}{\Delta x} \int_{x_{i-\frac12}}^{x_{i+\frac12}} \dd x\, s(q(t, x)) \label{eq:sourceexact}
\end{align}
as is usual in finite volume methods. For details on this, see e.g.\ \cite{barsukow18activeflux}. As Active Flux has access to point values $q_{i+\frac12}$ at cell interfaces, the integrals in \eqref{eq:fluxexact}--\eqref{eq:sourceexact} can immediately be evaluated via quadrature and there is no Riemann Problem to solve. In order to achieve third order in time, the quadrature in time needs values $q^n_{i+\frac12}, q^{n+\frac12}_{i+\frac12}, q^{n+1}_{i+\frac12}$ and then e.g. the flux is given by Simpson's rule
\begin{align}
 \hat f_{i+\frac12} = \frac{1}{6} \Big(f(q^n_{i+\frac12}) + 4 f(q^{n+\frac12}_{i+\frac12}) + f( q^{n+1}_{i+\frac12}) \Big) \label{eq:fluxquadrature}
\end{align}

In order to obtain the point values $q^{n+\frac12}_{i+\frac12}$ and $ q^{n+1}_{i+\frac12}$, a globally continuous reconstruction $\{q_{\text{recon},i}\}_i$, $q_{\text{recon},i} : [-\frac{\Delta x}{2},\frac{\Delta x}{2}] \to \admis$ (see Figure \ref{fig:af_states}) serves as initial data (described in Section \ref{ssec:usualrecon} below, as well as in e.g. \cite{vanleer77,eymann13,barsukow18activeflux}). The point values are computed via an exact or approximate evolution operator, i.e.\ by solving an initial value problem at the location of the point value\footnote{The point values $q^{n+\frac12}_{i+\frac12}, q^{n+1}_{i+\frac12}$ are obtained from the \emph{same} initial data (reconstruction at time $t^n$).}. For the shallow water equations, the (approximate) evolution operator is discussed in Section \ref{ssec:approxoperator}.
The full method consists of the update procedure for the point values, and \eqref{eq:generalafavgupdate} with \eqref{eq:fluxquadrature}.

\subsection{Usual reconstruction and limiting} \label{ssec:usualrecon}

Every cell has access to three pieces of information, and thus in every cell a parabola $q_{\text{recon},i}(x)$ ($x \in \left[-\frac{\Delta x}{2}, \frac{\Delta x}{2}\right]$) can be constructed which passes through the given point values $q_{\text{recon},i}\left(\pm\frac{\Delta x}{2}\right) =  q_{i\pm\frac12}$ and whose average agrees with the given cell average 
\begin{align}
\frac{1}{\Delta x} \int_{-\Delta x/2}^{\Delta x/2} \dd x \, q_{\text{recon},i}(x) = \bar q_i
\end{align}

\begin{align} \begin{split}
 q_{\text{recon},i}(x) &= -3 \left(2 \bar q_i - q_{i-\frac12} - q_{i+\frac12}\right) \frac{x^2}{\Delta x^2} \\&\hspace{-1cm}+ \left(q_{i+\frac12} - q_{i-\frac12}\right) \frac{x}{\Delta x} + \frac{6 \bar q_i - q_{i-\frac12} - q_{i+\frac12}}{4} \quad \text{for } x \in \left[-\frac{\Delta x}{2}, \frac{\Delta x}{2}\right] \end{split} \label{eq:parabolicrecon} 
\end{align} 

The individual parabolae meet at the cell interfaces, but their derivatives there generally do not agree. The parabolic reconstruction is one of the ingredients to achieve a third-order method. 

In order to prevent spurious oscillations (e.g.\ at discontinuities), limiting needs to applied. Suggestions in the literature include replacing the parabolic reconstruction in the cell by several parabolae or using a hyperbolic reconstruction instead (\cite{roe15,kerkmann18}). Additionally, discontinuous reconstructions have been considered in \cite{eymann11,eymann13a,kerkmann18} as limiting strategies, but they violate the philosophy of Active Flux and might spoil structure preservation. In \cite{barsukow19activeflux} a limiting strategy was suggested which replaces the parabolae by power laws. For what follows, it is convenient to define the reconstruction with respect to any interval $[x_\text{L}, x_\text{R}]$. In this paper, the following limited reconstruction is used:
\begin{align}
\mathscr R\left(\bar q, q_\text L, q_\text R; x_\text L, x_\text R; x\right) &= 
\begin{cases}
 \mathscr R^\text{limited}\left(\bar q,  q_\text L,q_\text R; x_\text L, x_\text R; x\right) & 
   \text{if } \left(\bar q,  q_\text L,q_\text R\right)\in \mathfrak L\\
 \mathscr R^\text{parabolic}\left(\bar q, q_\text L, q_\text R; x_\text L, x_\text R; x\right) & \text{else}
\end{cases} \label{eq:limitingpowerlaw}
\end{align}
where the limited reconstruction is applied if the states are in

\begin{multline}
\mathfrak L := \bigg\{
			\left(\bar q,  q_\text L,q_\text R\right)\in \mathbb R^3\,:\,
			\left(q_\text L < \bar q < q_\text L + \frac{q_\text R - q_\text L}{3}\right) 
			\\ \text{ or } \left(q_\text R - \frac{q_\text R - q_\text L}{3} < \bar q < q_\text R\right) 
			 \text{or } \left(q_\text R < \bar q < q_\text R + \frac{q_\text L - q_\text R}{3}\right) \\
			\text{ or } \left(q_\text L - \frac{q_\text L - q_\text R}{3} < \bar q < q_\text L\right) 
			\bigg\}
\end{multline}
The limited and parabolic reconstruction are given by
\begin{align}
&\mathscr R^\text{limited}\left(\bar q, q_\text L, q_\text R; x_\text L, x_\text R; x\right) = \\
&\phantom{mmmmm} \nonumber \begin{cases}
 q_\text L + (q_\text R - q_\text L) \left\lvert \displaystyle  \frac{x - x_\text L}{x_\text R - x_\text L} \right \rvert^{\mathscr E} & \text{if }\frac{1}{\mathscr E_\text{max}} \leq \mathscr E \leq 
 \mathscr E_\text{max}\\
 \mathscr R^\text{parabolic}\left(\bar q, q_\text L, q_\text R; x_\text L, x_\text R; x\right) & \text{else}
\end{cases}\\
\mathscr E &:= \frac{q_\text R - \bar q}{\bar q - q_\text L}
\end{align}
and
\begin{align}
\mathscr R^\text{parabolic}&\left(\bar q, q_\text L, q_\text R; x_\text L, x_\text R; x\right) = 
\frac{ q_\text R x_\text L^2 + 2 x_\text L x_\text R (-3 \bar q + q_\text R + q_\text L) + q_\text L x_\text R^2 }{(x_\text R-x_\text L)^2} 
\\ \nonumber &+
x \frac{ 6 \bar q (x_\text R+x_\text L) - 2 \Big(q_\text R(2 x_\text L + x_\text R) + q_\text L(2 x_\text R + x_\text L)\Big) }{ (x_\text R-x_\text L)^2 }
\\ \nonumber &+
3 x^2 \frac{ q_\text L -2 \bar q + q_\text R}{(x_\text R-x_\text L)^2}
\end{align}
This reconstruction fulfills
\begin{align}
 \mathscr R\left(\bar q, q_\text L, q_\text R; x_\text L, x_\text R; x_\text L\right) &= q_\text L \\
 \mathscr R\left(\bar q, q_\text L, q_\text R; x_\text L, x_\text R; x_\text R\right) &= q_\text R \\
 \int_{x_\text L}^{x_\text R} \dd x \, \mathscr R\left(\bar q, q_\text L, q_\text R; x_\text L, x_\text R; x\right) &= (x_\text R-x_\text L) \bar q \label{eq:reconconservation}
\end{align}
such that the usual reconstruction is $q^n_{\text{recon},i}(x) = \mathscr R\left(\bar q^n_i, q^n_{i-\frac12}, q^n_{i+\frac12}; -\frac{\Delta x}{2}, \frac{\Delta x}{2}; x\right)$.

We choose $\mathscr E_\text{max} = 50$ to avoid {excessive amplification of} rounding errors. For more details on this limiting procedure and also plots of the function for different combinations of $\bar q, q_\text L, q_\text R$, see \cite{barsukow19activeflux}.

\section{Bottom topography and non-negative reconstruction} \label{sec:bottomandrecon}

\subsection{Overview of design decisions}

Different ingredients need to be brought together in order to endow the method with all the following properties:

\begin{itemize}
 \item third order accurate
 \item well-balanced on wet states
 \item positivity preserving
 \item allowing for dry states
 \item well-balanced on partially wet states
\end{itemize}

Certain design decisions turn out not to be independent from each other. In this section therefore, a general reasoning is given. The actual formulae are given in the sections below, and certain technical details are even deferred to the appendix.

Third order of accuracy has been so far achieved for Active Flux by considering a piecewise parabolic reconstruction of the conserved quantities. For the application to shallow water equations, the method needs to keep the lake at rest stationary. We decide to \textbf{project the bottom topography} onto a globally continuous, piecewise parabolic function. This ensures that the lake at rest is reconstructed exactly. An exact evolution operator would then keep this state stationary, and one might hope to find an approximate evolution operator that is still exact on the lake-at-rest state. This line of thought is pursued in this paper. However, it is incomplete yet, and needs to be augmented by further ingredients necessary for hyperbolic problems.

A parabolic reconstruction requires a limiting procedure. However, this might entail limiting of the (locally) parabolic water height $h$ of a lake at rest! In this paper, we thus decide to \textbf{reconstruct the equilibrium variable $h+b$} instead, and to subtract the bottom from it to obtain a reconstruction of $h$. In some cell, when $h+b$ is reconstructed parabolically and the parabolic $b$ is subtracted from the result, the reconstruction of $h$ is just the same, as if $h$ would have been directly reconstructed parabolically. The only difference is when limiting turns on. 

Having brought well-balancing and limiting into concordance, non-negativity pre\-servation needs to be included. The usual limiting is not sufficient. We take the decision not to aim at detecting violations of non-negativity precisely, but to be slightly over-careful. This makes the enforcement of \textbf{non-negativity independent of the choice of limiting}. Then, in cells where a possible violation of non-negativity is detected, $h$ (and not $h+b$) is reconstructed in a particular fashion described later. It is chosen to ensure that half-wet cells at the shore of a lake at rest are also reconstructed exactly.

\subsection{Bottom topography}

We replace the exact bottom topography $b(x)$ by a piecewise parabolic approximant $\hat b(x)$:
\begin{align}
 \hat b(x) &= b(x_{i-\frac12}) 
 + \left(x - \frac{x_{i+\frac12} + x_{i-\frac12}}{2}\right) \frac{ b(x_{i+\frac12}) - b(x_{i-\frac12})}{x_{i+\frac12} - x_{i-\frac12}} \nonumber
 \\&\phantom{mmmm}+ 2 \left(x - \frac{x_{i+\frac12} + x_{i-\frac12}}{2} \right)^2 \frac{ b(x_{i+\frac12}) - 2 b(x_i) + b(x_{i-\frac12}) }{\left(x_{i+\frac12} - x_{i-\frac12}\right)^2} \\
 \nonumber &\phantom{mmmm} \qquad \text{if } x_{i-\frac12} \leq x < x_{i+\frac12}
\end{align}

Only the projected bottom is used in the following, and the notational difference is suppressed.

The projection onto globally continuous, piecewise quadratic polynomials is not strictly necessary, but simplifies certain aspects and the imple\-men\-ta\-tion of the method. A discontinuous bottom topography should be regularized prior to usage.

\subsection{Non-negative reconstruction}

\begin{figure}
 \centering
 \includegraphics[width=\textwidth]{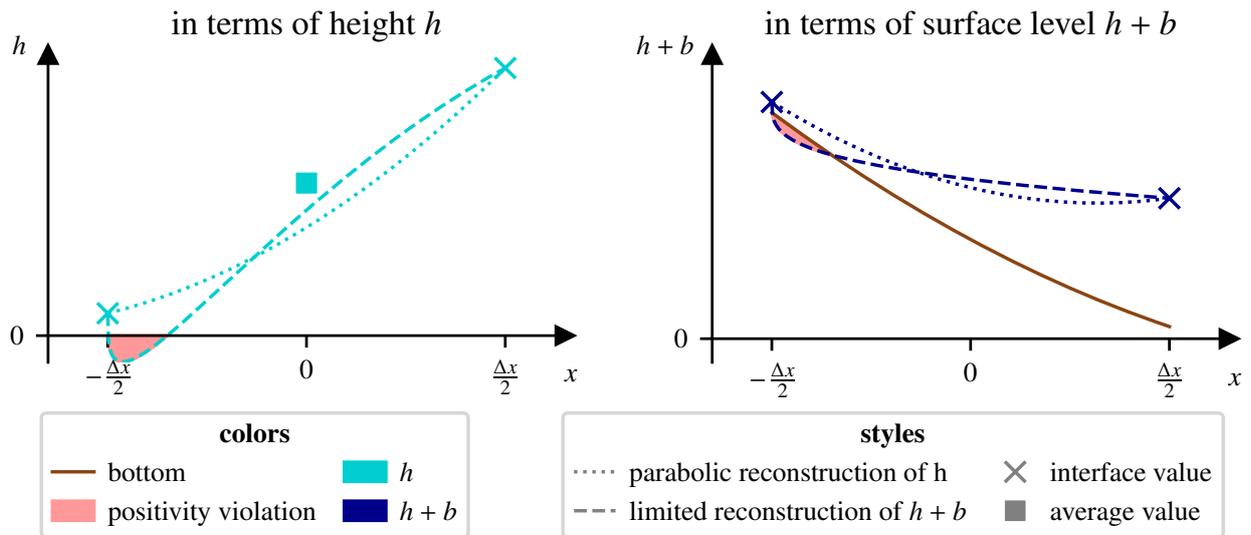}
 \caption{As the bottom topography (\emph{brown}) is parabolic, reconstructing $h+b$ parabolically and subtracting $b$ gives the same result as reconstructing $h$ parabolically straight away (\emph{light blue}). A difference arises only when limiting is involved, i.e.\ when the reconstruction is non-parabolic. In the example shown here, the parabolic reconstruction of $h+b$ is positive, but non-monotone and requires limiting. The limited reconstruction of $h+b$ (\emph{dark blue}), however, violates non-negativity. Note that all the reconstruction procedures are conservative (the average water height is indicated by the square symbol). \emph{Left}: The surface levels resulting from the different reconstructions are shown. \emph{Right}: The reconstructed water height is shown. The two figures differ only by adding or subtracting $b(x)$ to the curves.}
 \label{fig:limitingpositivityviolation}
\end{figure}

\begin{figure}
  \centering
  \includegraphics[width=\textwidth]{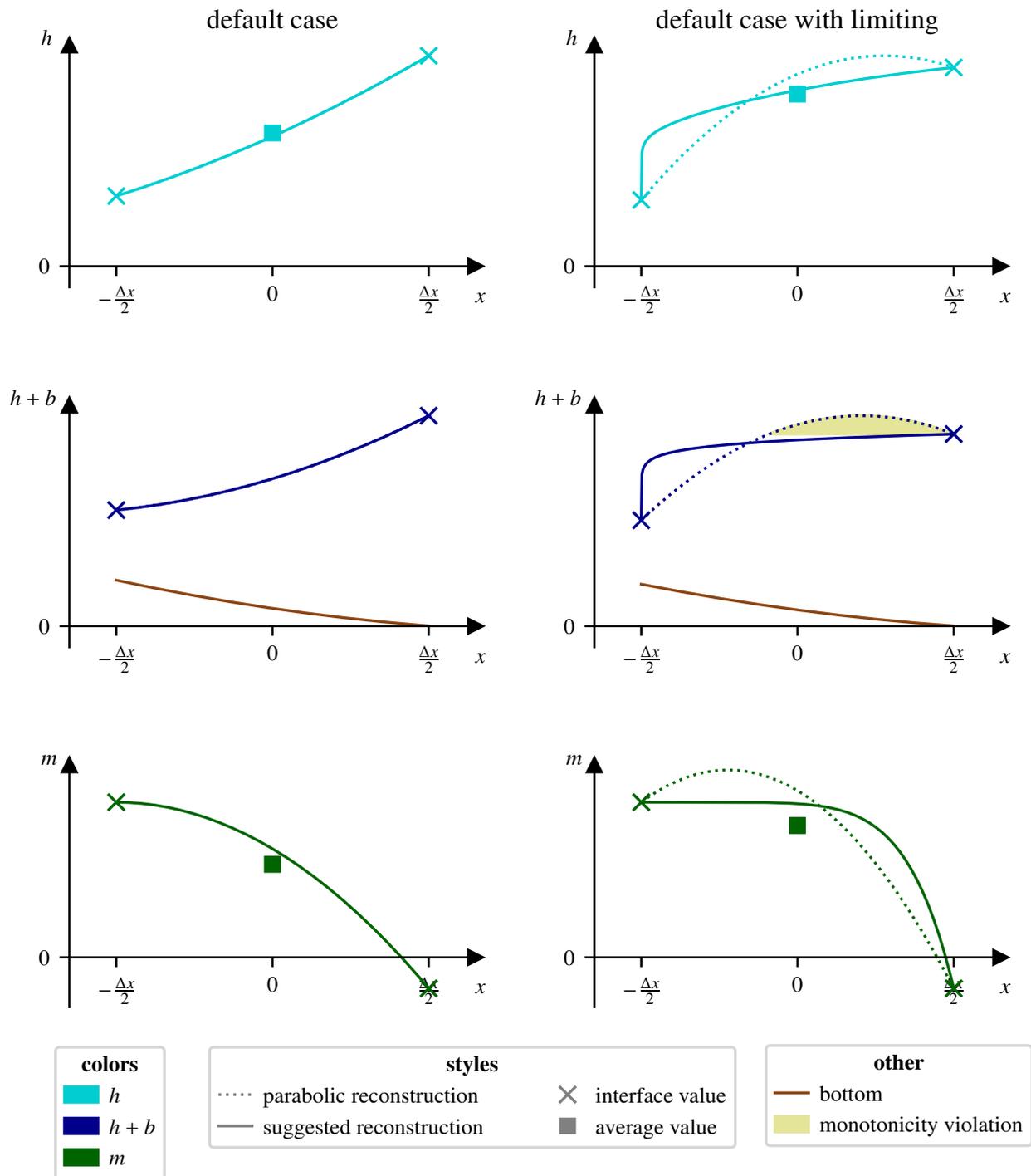}
  \caption{Usual reconstruction, used when the water height remains positive. \emph{Left}: Parabolic reconstruction. \emph{Right}: Limited reconstruction. \emph{Top}: Water height. \emph{Center}: Water level. \emph{Bottom}: Momentum.}
  \label{fig:casedefault}
 \end{figure}

The reconstruction of the water height needs to remain non-negative -- assuming that the point values and the cell average are non-negative. First, non-negativity of a parabolic reconstruction of $h$ is checked. If a parabolic reconstruction of $h$ would violate non-negativity of $h$, then $h$ and $m$ are reconstructed in a particular fashion as described below. If the parabolic reconstruction of $h$ is positive, we study whether a limited reconstruction of $h+b$ would violate non-negativity of $h$ (see Figure \ref{fig:limitingpositivityviolation}). To make sure that limiting of $h+b$ does not spoil non-negativity of $h$ we suggest the following condition:

If 
\begin{align}
 \min(h_\text L + b_\text L, h_\text R + b_\text R) > \max(b_\text L, b_\text R)
\end{align}
then reconstruct (see Figure \ref{fig:casedefault})
\begin{align}
  h_\text{recon}(x) &= \mathscr R\left(\bar h + \bar b, h_\text L + b_\text L, h_\text R + b_\text R; -\frac{\Delta x}{2}, \frac{\Delta x}{2}; x\right) - b(x)
 \end{align}
with $\bar b = \frac{1}{6} (b_\text L + 4 b(\frac{x_\text L + x_\text R}{2}) + b_\text R)$; otherwise reconstruct $h$ directly:
 \begin{align}
  h_\text{recon}(x) &= \mathscr R\left(\bar h, h_\text L, h_\text R; -\frac{\Delta x}{2}, \frac{\Delta x}{2}; x\right)
 \end{align}
In both these cases,
 \begin{align}
  m_\text{recon}(x) &= \mathscr R\left(\bar m, m_\text L, m_\text R; -\frac{\Delta x}{2}, \frac{\Delta x}{2}; x\right)
 \end{align}
Appendix \ref{app:reconstruction} contains a flow chart of the reconstruction procedure.

The following lemma establishes when a parabolic reconstruction has negative values in some region inside the cell:

\begin{lemma} \label{lem:nonnegativityparabolic}
Assume $h_\text L > 0$, $h_\text R > 0$, $\bar h > 0$. There exists a region inside $[-\Delta x/2, \Delta x/2]$ where the parabolic reconstruction \eqref{eq:parabolicrecon} becomes negative if and only if
\begin{align}
 3 \bar h < h_\text L + h_\text R - \sqrt{h_\text Lh_\text R} \quad \text{and}\quad \bar h < \frac{h_\text L + h_\text R}{2} - \frac{ \lvert h_\text R - h_\text L \rvert }{6 }  \label{eq:reconnonneg}
\end{align}
simultaneously.
\end{lemma}
\begin{proof}
\begin{enumerate}[i)]
\item If the parabolic reconstruction becomes negative inside the cell, then this must happen in the interior as $h_{\text L/\text R} > 0$. The (unique) extremum of \eqref{eq:parabolicrecon} is located at
\begin{align}
 \Delta x \frac{ h_\text R - h_\text L }{6 (2 \bar h - h_\text L - h_\text R)}
\end{align}
where the reconstruction attains the value
\begin{align}
  \frac{ (h_\text R - h_\text L)^2 }{12 (2 \bar h - h_\text L - h_\text R)}  + \frac{6 \bar h - h_\text L - h_\text R}{4} 
\end{align}
The extremum is located inside $[-\Delta x/2, \Delta x/2]$ iff
\begin{align}
 \frac{h_\text L + h_\text R}{2} - \frac{ \lvert h_\text R - h_\text L \rvert }{6 } \geq \bar h 
\end{align}
where $2 \bar h - h_\text L - h_\text R < 0$ has been used, as otherwise the reconstruction is concave and the extremum is a maximum.

This value being less than zero implies
\begin{align}
  \frac{ (h_\text R - h_\text L)^2 }{12 }  + \frac{\left(2 \bar h - (h_\text L + h_\text R)\right)(6 \bar h - (h_\text L + h_\text R))}{4} > 0\\
  h_\text R^2  + h_\text R h_\text L +  h_\text L^2   + \left(3\bar h\right) ^2 - 2 \left(3\bar h\right) (h_\text L + h_\text R)   > 0
\end{align}
Thus
\begin{align}
 \left(3\bar h\right)  \lessgtr  (h_\text L + h_\text R) \mp \sqrt{  h_\text R h_\text L  }
\end{align}
Only the upper choice of signs is relevant.

\item Conversely, if \eqref{eq:reconnonneg} is true, then by the inequality of the geometric and arithmetic means
\begin{align}
 \bar h < \frac{h_\text L + h_\text R}{6} < \frac{h_\text L + h_\text R}{2}
\end{align}
such that the function is convex and the rest follows from the above. This completes the proof.
\end{enumerate}
\end{proof}

{\emph{Note}: As is discussed below, the case, in which $h_\text L = 0$, or $h_\text R = 0$, or $\bar h = 0$ will not be attempted to be reconstructed parabolically anyway.}

If the parabolic reconstruction turns negative inside the cell, we ignore the possibility that regular limiting might ``accidentally'' correct this. Instead, we choose an entirely different reconstruction. Given the point values $h_\text R$, $h_\text L$ and the average $\bar h$, the suggested non-negative reconstruction of the water height $h$ is detailed below. In principle, the reconstruction of the momentum $m$ is not subject to any a priori bounds. However, it has been found in practice that enforcing the positivity of the water height while not restricting in any way the reconstruction of the momentum can lead to unexpectedly high water speeds. Therefore we suggest to also adjust the reconstruction of the momentum to the reconstruction of the water height.

\newcommand{\linto}{\overset{\text{lin}}{\to}}
\begin{definition}
Define the shorthand notation
 \begin{align}
  f(x) = a \linto b \quad x_1 \leq x \leq x_2
 \end{align}
 to denote a linear function $f : [x_1,x_2] \to [a,b]$ which takes the value $a$ at $x_1$ and $b$ at $x_2$, i.e.\ the function
 \begin{align}
  f(x) = a + (x - x_1) \frac{b - a}{x_2-x_1}
 \end{align}
\end{definition}
The non-negative reconstruction is computed as follows:

\begin{enumerate}[a.]

 \item \label{it:lowerthanboth} If $\bar h < h_\text R$ and $\bar h < h_\text L$ and if the parabolic reconstruction is violating non-negativity (see \eqref{eq:reconnonneg}), then we suggest to reconstruct $h$ by three linear segments joined in a continuous fashion, with the central segment being horizontal and having the value $f \bar h$ with $f \in (0,1)$, given below. 
 
 \begin{figure}
  \centering
  \includegraphics[width=\textwidth]{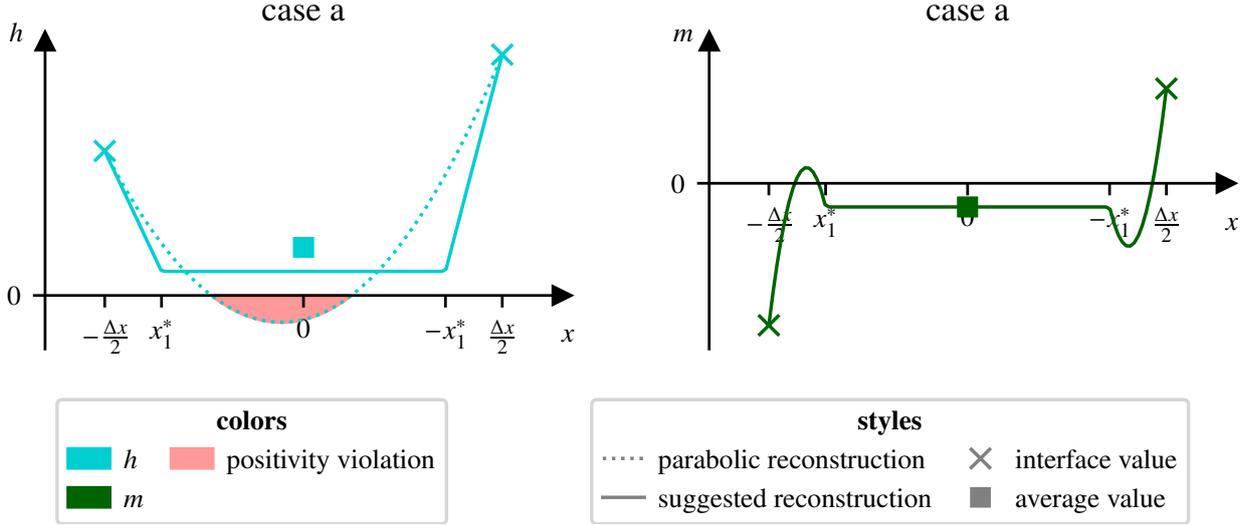}
  \caption{Reconstruction of the case $\bar h < h_\text R$ and $\bar h < h_\text L$, when a parabolic reconstruction would violate non-negativity. \emph{Left}: Water height. \emph{Right}: Momentum.}
  \label{fig:moses}
 \end{figure}

 Compute
 \begin{align}
    x^*_1 := \Delta x  \frac{ 2 \bar h   (f-2) + h_\text L + h_\text R}{4 \bar h f - 2 (h_\text L + h_\text R)} 
    \label{eq:nonnegavgbelowbothxlocation}
 \end{align}
 
 Then the following reconstruction
 \begin{align}
   h_\text{recon}(x) = 
\begin{cases} 
 h_\text L \linto f\cdot \bar h &   -\frac{\Delta x}{2} < x \leq x^*_1\\
 f\cdot \bar h & x^*_1 < x \leq -x^*_1 \\
 f\cdot \bar h \linto h_\text R & -x^*_1 < x \leq \frac{\Delta x}{2}
\end{cases}
\end{align}
is conservative.

Lemma \ref{lem:positivewaterheight}\ref{it:lemx1} (in the Appendix) demonstrates that $- \frac{\Delta x}{2} <  x^*_1 \leq 0$ and that this reconstruction therefore is well-behaved if
\begin{align}
 f \geq 2 - \frac{h_\text L + h_\text R}{2 \bar h}
\end{align}

We use $f = \max\left(\frac12,2 - \frac{h_\text L + h_\text R}{2 \bar h} \right) $ in practice. 

We suggest to reconstruct the momentum such that in the interval $[x^*_1, -x^*_1]$ the momentum is constant as well and we choose this constant to be the cell average $\bar m$:

\begin{align} m_\text{recon}(x) = \begin{cases}
      \mathscr R\left(\bar m, m_\text L, \bar m; -\frac{\Delta x}{2}, x^*_1; x\right) & -\frac{\Delta x}{2} < x \leq x^*_1\\
     \bar m & x^*_1 < x \leq  -x^*_1\\
     \mathscr R\left(\bar m, \bar m, m_\text R;  -x^*_1, \frac{\Delta x}{2}; x\right) & -x^*_1 < x \leq \frac{\Delta x}{2}
   \end{cases} \label{eq:mosesmomentum} \end{align}
   By Lemma \ref{lem:positivewaterheight}\ref{it:lemmosesmom} this reconstruction is conservative (see Fig. \ref{fig:moses}).

\item \label{it:heigtRlessH} 

If $h_\text R \leq \bar h < h_\text L$, then action is needed if the reconstruction is violating non-negativity (see \eqref{eq:reconnonneg}) or if $h_\text R = 0$. The latter case ($h_\text R = 0$) needs to be considered in order to reconstruct the lake at rest exactly, even if the parabolic reconstruction might be well-behaved in this case. 

In both cases we suggest to reconstruct the water height in a piecewise fashion -- taking the water height to be constant and equal to $h_\text R$ in the interval $( x^*,\frac{\Delta x}{2})$ (for some $x^*$ to be specified), and to reconstruct the remaining interval $(-\frac{\Delta x}{2}, x^*)$ such that $h+b$ is linear in it:

\begin{align}
 h_\text{recon}(x) &= \begin{cases}
                      h_\text L + b_\text L + \left(x + \frac{\Delta x}{2}\right) \frac{h_\text R + b(x^*) - h_\text L - b_\text L}{x^* + \frac{\Delta x}{2}} - b(x) & -\frac{\Delta x}{2} \leq x < x^*\\
                      h_\text R & x^* \leq x < \frac{\Delta x}{2}
                     \end{cases} 
\end{align}

In Lemma \ref{lem:recondrywetmirror} (in the Appendix), under certain conditions the existence of a suitable $x^*$ that makes this reconstruction conservative is established. For example, if the data belong to a lake at rest, then clearly there exists such an $x^*$. The equation governing $x^*$ is an algebraic equation of third degree and is given in the proof of the Lemma (see Figure \ref{fig:casebrightsmall}).

In a certain case (involving a concave bottom), a suitable $x^*$ does not exist, while the parabolic reconstruction would violate non-negativity and cannot be used. For this case, the above-mentioned Lemma gives an alternative reconstruction, using a water height linear (instead of constant) in space on the left side. For details and the computation of $x^*$, the reader is referred to Lemmas \ref{lem:recondrywet}--\ref{lem:recondrywetmirror} in the Appendix and Figures \ref{fig:casebrightsmallconcave} and \ref{fig:recondrywetLeft}. 

 \begin{figure}[p]
  \centering
  \includegraphics[width=0.95\textwidth]{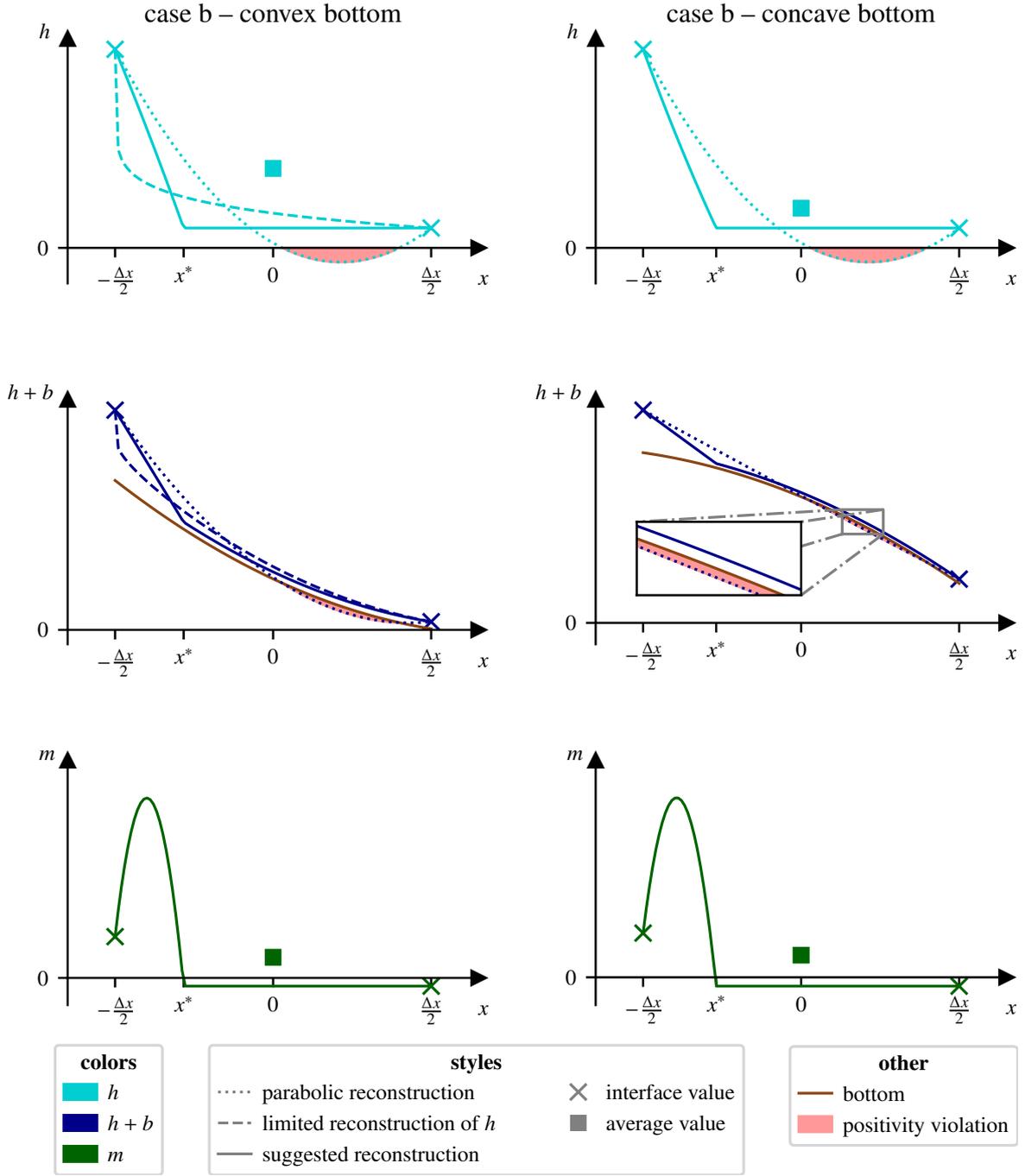}
  \caption{The non-negative reconstruction procedure for the case $h_\text R < \bar h < h_\text L$. \emph{Left}: Convex bottom. \emph{Right}: The case of a concave bottom in certain cases requires a modification to this reconstruction (see Figure \ref{fig:casebrightsmallconcave}); here the standard case is depicted. \emph{Top}: The water height $h$, consisting of a constant and a parabolic piece, joined continuously. The parabolic piece is such that $h+b$ is linear. \emph{Center}: The water surface $h+b$. \emph{Bottom}: Associated reconstruction of the momentum. Several aspects of the reconstruction are exemplified in the figure. First, the top figure on the left demonstrates that the limiter sometimes can handle a positivity violation of the parabola. In order to decouple limiting from positivity, we decide never to rely on this and to replace the parabola by the piecewise defined positive reconstruction even if the limiter would yield a positive reconstruction. Second, the reconstruction of the momentum is chosen constant in the interval where the reconstruction of $h$ is chosen constant.}
  \label{fig:casebrightsmall}
 \end{figure}

 \begin{figure}
  \centering
  \includegraphics[width=0.95\textwidth]{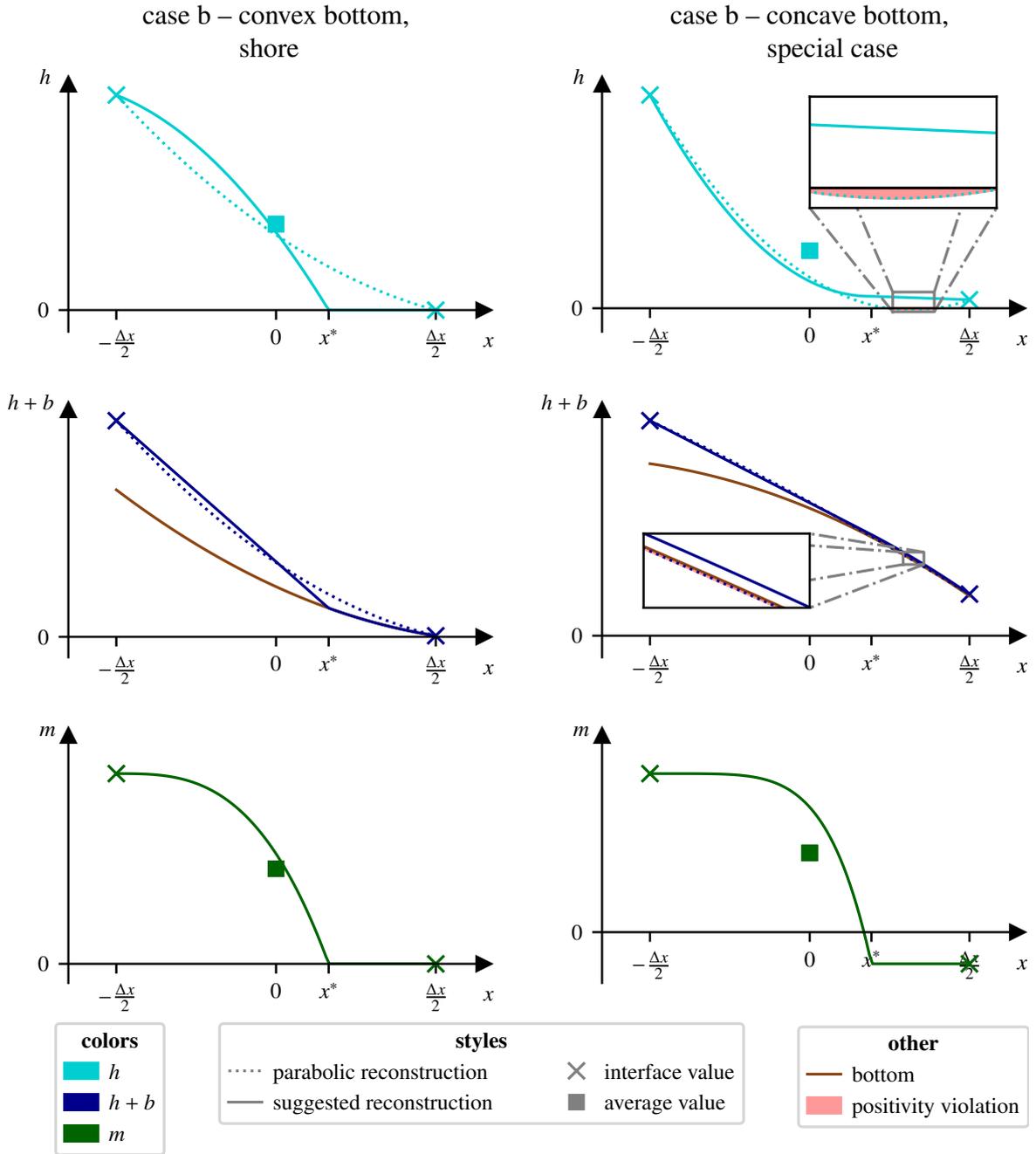}
  \caption{Special cases of the non-negative reconstruction procedure for the case $h_\text R < \bar h < h_\text L$. \emph{Top}: The water height $h$. \emph{Center}: The water surface $h+b$. \emph{Bottom}: Associated reconstruction of the momentum. \emph{Left}: $h_\text R = 0$. Because $h$ is reconstructed using a constant piece and a parabolic piece that turns $h+b$ linear, and because the momentum is chosen constant where $h$ is reconstructed constant, this case does not require any special treatment. \emph{Right}: The exceptional case when the bottom is concave. Unfortunately, it occurs in a very narrow range of parameters which makes it difficult to visualize. Observe that the piece of $h_\text{recon}$ that would usually be chosen constant, now is linear with a non-vanishing slope. For more details on this case see Appendix \ref{app:reconstruction} and also Figure \ref{fig:recondrywetLeft}.}
  \label{fig:casebrightsmallconcave}
 \end{figure}

We suggest the following momentum reconstruction:

   \begin{align} m_\text{recon}(x) = \begin{cases}
	\mathscr R\left(\frac{ \bar m  - \left(\frac{1}{2} - \frac{x^*}{\Delta x}\right) m_\text{R}}{\frac{1}{2} + \frac{x^*}{\Delta x}}, m_\text{L}, m_\text R; -\frac{\Delta x}{2}, x^*; x\right) & -\frac{\Delta x}{2}\leq x < x^*\\
	m_\text R & x^* \leq x < \frac{\Delta x}{2} 
   \end{cases} \label{eq:recononeflatx2} \end{align}

By Lemma \ref{lem:positivewaterheight}\ref{it:lemhalf1} this reconstruction is conservative. Note that the evolution operator described below guarantees that $h_\text R = 0$ implies $m_\text R = 0$.
   
 \item If $h_\text L \leq \bar h < h_\text R$ and either $h_\text L = 0$ or if the parabolic reconstruction is not non-negative, \eqref{eq:reconnonneg}, then the situation is analogous and a mirror image of the one discussed previously. We suggest to reconstruct the water height as

 \begin{align}
  h_\text{recon}(x) = \begin{cases} h_\text L & -\frac{\Delta x}{2} \leq x < x^*\\
                       h_\text L + b(x^*) + (x - x^*) \frac{h_\text R + b_\text R - h_\text L - b(x^*)}{\frac{\Delta x}{2} - x^*} - b(x) & x^* \leq x < \frac{\Delta x}{2}
                      \end{cases}
 \end{align}
 with an exceptional case discussed again in the Appendix (Lemma \ref{lem:recondrywet}).

The momentum reconstruction reads:
  \begin{align} m_\text{recon}(x) = \begin{cases}
  m_\text L & -\frac{\Delta x}{2} \leq x < x^* \\
  \mathscr R\left(\frac{\bar m  -  m_\text L \left (\frac{x^*}{\Delta x} + \frac{1}{2}\right)}{\frac{1}{2} - \frac{x^*}{\Delta x}}, m_\text L, m_\text R; x^*, \frac{\Delta x}{2}; x\right) & x^* \leq x < \frac{\Delta x}{2}
  \end{cases} \label{eq:recononeflatx3} \end{align}

\end{enumerate}
This completes the presentation of the non-negative reconstruction.
If none of the above cases applies, then it is safe to reconstruct $h$ parabolically (or, equivalently, to reconstruct $h+b$ parabolically and to subtract $b$). 
A flow chart of the reconstruction procedure is given in Appendix \ref{app:reconstruction}. As an example, Figure \ref{fig:doublerarefaction} shows which cells are reconstructed how in a representative simulation containing dry areas, shocks and rarefactions.

\section{Update of point values} \label{sec:pointvalues}
\subsection{Overview}

Given a non-negative reconstruction, the point values are evolved using an approximate evolution operator based on those introduced in \cite{barsukow19activeflux} and \cite{barsukow19activefluxsource}. This operator is described in Section \ref{ssec:approxoperator}. 

For practical usage it requires three modifications. As noticed in \cite{kerkmann18,barsukow19activeflux} standard choices of evolution operators for Active Flux can lead to artificially stationary shocks (which violate the Rankine-Hugoniot conditions) and to artefacts in transonic rarefactions, such as those well-known for finite volume methods. In \cite{barsukow19activeflux}, an \textbf{entropy fix} is proposed which removes both artefacts; the way it is applied here to the shallow water equations is described in Section \ref{ssec:entropyfix}. This is the first modification. Modifications related to the presence of \textbf{dry, or half-wet cells} are discussed in Section \ref{ssec:nearvacuum} and Section \ref{ssec:drycellsregularizationpointvalues}. The third modification guarantees that the lake-at-rest state is resolved to machine precision (\textbf{well-balancing}). This is described in Section \ref{ssec:wellbalancing}.

\subsection{Approximate evolution operator} \label{ssec:approxoperator}

As outlined in \cite{barsukow19activeflux}, in order for Active Flux to achieve third order, a sufficiently accurate approximate solution $\tilde q(t,x)$ is required, which approximates the exact solution $q(t, x)$ for every $x$ as follows
\begin{align}
 \tilde q(t, x) = q(t, x) + \mathcal O(t^3)
\end{align}

Consider the following system 
\begin{align}
 \del_t Q_i + \lambda_i(Q_1, \ldots, Q_M)  \del_x Q_i = S_i(x; Q_1, \ldots, Q_M) \qquad i = 1, \ldots, M \label{eq:diagonalwithsource}
\end{align}
which arises upon diagonalization of the Jacobian in \eqref{eq:hypsystemgeneral}, if it is possible to find characteristic variables {$Q$. In the following, the existence of such a transformation from conservative variables $q$ to the characteristic variables $Q$ is assumed, which is the case for the shallow water equations. For general hyperbolic systems, an analogous approximate evolution operator has been derived in \cite{barsukow19activeflux}.} The corresponding initial data are denoted by $Q_{i,0}(x)$, $i = 1, \ldots, M$.

\begin{example}
For the shallow water equations, $M = 2$ and
\begin{align}
 Q_\pm &= 2c \pm v &  \lambda_\pm &= v \pm c & S_\pm &= \mp g \del_x b \label{eq:charvarshallowwater}
\end{align}
with $c = \sqrt{g h}$, {having renamed $Q_1 \equiv Q_+$ and $Q_2 \equiv Q_-$, because for this specific system, $Q_\pm$ are more telling names (and analogously $\lambda_\pm$).}
\end{example}

Local linearization, i.e.
\begin{align}
 \frac{ Q_i(t, x) - Q_{i,0}(x - \lambda_i(x)) }{ t} = S_i(x)
\end{align}
with the shorthand notation
\begin{align}
 {S_i} \equiv S_i(x; Q_{1,0}(x), \ldots, Q_{M,0}(x)) \label{eq:evoopshorthandx}\\
 {\lambda_i} \equiv \lambda_i(Q_{1,0}(x), \ldots, Q_{M,0}(x)) \label{eq:evoopshorthandx2}
\end{align}
yields a solution operator that is one order less accurate than required. Consider therefore
\begin{align}
 \frac{ Q_i(t, x) - Q_{i,0}(x - \lambda_i^* t) }{t} = S_i^*
\end{align}
instead, with
\begin{align}
 \lambda_i^* &:= \lambda_i\left(Q_1\left(\frac{t}{2}, x - \lambda_i \frac{t}{2}\right), \ldots, Q_M\left(\frac{t}{2}, x - \lambda_i \frac{t}{2} \right )\right ) + \mathcal O(t^2) \label{eq:approxevolambdaprelim}\\
 S_i^* &:= S_i\left(x - \lambda_i \frac{t}{2}; Q_1\left(\frac{t}{2}, x - \lambda_i \frac{t}{2}\right), \ldots, Q_M\left(\frac{t}{2}, x - \lambda_i \frac{t}{2}\right)\right) + \mathcal O(t^2) \label{eq:approxevosourceprelim}
\end{align}
This yields the desired order of accuracy:

\begin{theorem} \label{thm:evolutionoperator}
 Consider the approximate solution operator
\begin{align}
 \tilde Q_i(t, x) = Q_{i,0}(x - \lambda_i^* t) + t S_i^* \label{eq:evolutionoperatorprelim}
\end{align}
with \eqref{eq:approxevolambdaprelim}--\eqref{eq:approxevosourceprelim}. Then $\tilde Q_i(t, x) = Q_i(t, x) + \mathcal O(t^3)$.
\end{theorem}
\textit{Proof}. 
 All summations involve indices from 1 to $M$. Compute first

  \begin{align}
  \lambda_i^*  \rvert_{t=0} &= \lambda_i \qquad S_i^* \rvert_{t=0} = S_i\\
  \del_t \lambda_i^*  \rvert_{t=0} &= \sum_\alpha \frac{\del \lambda_i}{\del Q_\alpha} \frac12 \left(\del_t Q_\alpha\Big\rvert_{t=0} \!\!\!\! - \lambda_i Q'_{\alpha,0}(x) \right)\\
  &= \sum_\alpha \frac{\del \lambda_i}{\del Q_\alpha} \frac12 \Big(- (\lambda_\alpha + \lambda_i)  Q'_{\alpha,0}(x)  + S_\alpha \Big)\\
  \del_t S_i^*  \rvert_{t=0} &=-\frac12 \lambda_i \del_x S_i + \sum_\alpha \frac{\del S_i}{\del Q_\alpha}  \frac12 \Big(- (\lambda_\alpha + \lambda_i) Q'_{\alpha,0}(x)  + S_\alpha \Big)
 \end{align}
 Also,
 \begin{align}
 \del_t \tilde Q_i(t, x) &= Q_{i,0}'(x - \lambda_i^* t) \Big(- (\del_t \lambda_i^*) t - \lambda_i^*\Big) + S_i^* + t \del_t S_i^*\\
  \left. \del_t \tilde Q_i(t, x) \right \rvert_{t=0} &=  -Q_{i,0}'(x)  \lambda_i + S_i
 \end{align}
 \begin{align}
  \left. \del_t^2 \tilde Q_i(t, x) \right \rvert_{t=0} &= 
  Q_{i,0}''(x) \lambda_i^2
  - 2 Q_{i,0}'(x ) \del_t \lambda_i^*  \rvert_{t=0}
  + 2 \del_t S_i^* \rvert_{t=0}\\
  &= Q_{i,0}''(x) \lambda_i^2
  -  Q_{i,0}'(x ) \sum_\alpha \frac{\del \lambda_i}{\del Q_\alpha}  \Big(- (\lambda_\alpha + \lambda_i)  Q'_{\alpha,0}(x)  + S_\alpha \Big)
  \\\nonumber&- \lambda_i \del_x S_i + \sum_\alpha \frac{\del S_i}{\del Q_\alpha}  \Big(- (\lambda_\alpha + \lambda_i)  Q'_{\alpha,0}(x)  + S_\alpha \Big)
 \end{align}
 
 On the other hand, from \eqref{eq:diagonalwithsource}, by performing the Cauchy-Kowalevskaya procedure
 \begin{align}
  \del_t Q_i &= - \lambda_i \del_x Q_i + S_i\\
  \del_t^2 Q_i\rvert_{t=0} &=\left. \left( - \sum_\alpha \frac{\del \lambda_i}{\del Q_\alpha} \del_t Q_\alpha \del_x Q_i - \lambda_i \del_x (\del_t Q_i) + \sum_\alpha \frac{\del S_i}{\del Q_\alpha} \del_t Q_\alpha \right )\right \rvert_{t=0}\\
  &= - \sum_\alpha \frac{\del \lambda_i}{\del Q_\alpha} (- \lambda_\alpha  Q'_{\alpha,0} + S_\alpha)  Q'_{i,0} - \lambda_i \del_x (- \lambda_i Q'_{i,0} + S_i) \\\nonumber& + \sum_\alpha \frac{\del S_i}{\del Q_\alpha} (- \lambda_\alpha Q'_{\alpha,0} + S_\alpha)\\
  &= - \sum_\alpha \frac{\del \lambda_i}{\del Q_\alpha} (- \lambda_\alpha  Q'_{\alpha,0} + S_\alpha)  Q'_{i,0} 
  + \sum_\alpha \frac{\del \lambda_i}{\del Q_\alpha} \lambda_i Q'_{\alpha,0} Q'_{i,0}  + \lambda_i^2 Q''_{i,0} \\\nonumber& - \lambda_i \del_x S_i - \lambda_i \sum_\alpha \frac{\del S_i}{\del Q_{\alpha}} Q'_{\alpha,0}
  + \sum_\alpha \frac{\del S_i}{\del Q_\alpha} (- \lambda_\alpha Q'_{\alpha,0} + S_\alpha)\\
  &=  \sum_\alpha \frac{\del \lambda_i}{\del Q_\alpha} \Big( (\lambda_\alpha  + \lambda_i ) Q'_{\alpha,0}  - S_\alpha \Big) Q'_{i,0}  + \lambda_i^2 Q''_{i,0} \\\nonumber&  + \sum_\alpha \frac{\del S_i}{\del Q_\alpha} \Big( - (\lambda_i + \lambda_\alpha)  Q'_{\alpha,0} + S_\alpha \Big) - \lambda_i \del_x S_i
 \end{align}

Equations \eqref{eq:approxevolambdaprelim}--\eqref{eq:approxevosourceprelim}, as they are given, involve terms that are not available. However, as they are only required up to $\mathcal O(t^2)$, one can use local linearization as a predictor step:
\begin{align}
 Q_j\left(\frac{t}{2}, x - \lambda_i \frac{t}{2}\right) &\simeq Q_{j,0}\left(x - \frac{\lambda_i + \lambda_j}{2} t \right) \nonumber
 \\&\phantom{mmm}+ \frac{\Delta t}{2} S_{j}\left(x; Q_{1,0}\left(x \right), \ldots, Q_{M,0}\left(x  \right)\right) 
\end{align}

{The final shape of the approximate evolution operator for the update of the point values, therefore, is
\begin{align}
 \tilde Q_i(t, x) = Q_{i,0}(x - \lambda_i^* t) + t S_i^* \label{eq:evolutionoperator}
\end{align}
with
\begin{align}
 \lambda_i^* &:= \lambda_i\left(Q^*_{1,i}, \ldots, Q^*_{M,i}\right ) \label{eq:approxevolambda}\\
 S_i^* &:= S_i\left(x - \lambda_i \frac{t}{2}; Q^*_{1,i}, \ldots, Q^*_{M,i}\right) \label{eq:approxevosource}
\end{align}
and
\begin{align}
 Q^*_{j,i} := Q_{j,0}\left(x - \frac{\lambda_i + \lambda_j}{2} t \right)
 + \frac{\Delta t}{2} S_{j}\Big(x; Q_{1,0}\left(x \right), \ldots, Q_{M,0}\left(x  \right)\Big) \label{eq:predictor}
\end{align}
$\lambda_i$ and $S_i$ are given by the shorthand notation of \eqref{eq:evoopshorthandx}--\eqref{eq:evoopshorthandx2}
}

In the homogeneous case, this operator agrees with one of the operators suggested in \cite{barsukow19activeflux}. For the case of a linear homogeneous part, the operator suggested in \cite{barsukow19activefluxsource} is similar to the one above, but uses a slightly different predictor step. This is because the operator in \cite{barsukow19activefluxsource} has been derived as a Runge-Kutta scheme. For details, the reader is referred to the respective works.

{Boundary conditions can be easily implemented e.g. through the usage of ghost cells.}

The usage of characteristics raises the question of whether this solution operator is suitable in presence of shocks. Indeed, in \cite{barsukow19activeflux} it has been shown that a modification is necessary to account for the possibility of crossing characteristics. Otherwise, non-entropic features arise. By analogy to similar concepts for finite volume schemes, this modification is referred to as entropy fix and discussed next.

\subsection{Entropy fix} \label{ssec:entropyfix}

The derivation of the approximate evolution operator as presented in Section \ref{ssec:approxoperator} assumes a smooth solution. For hyperbolic systems of conservation laws, the solution is known to develop discontinuities in finite time. This is because characteristics can cross. The natural modification, suggested in \cite{barsukow19activeflux} for scalar problems initially, is to consider several candidate characteristics and use one of them based on some selection criterion. It does not seem necessary to compute exactly which of the characteristics to use. Here, in parallel to what has been suggested in \cite{barsukow19activeflux}, two candidates are computed by replacing $\lambda_i(x)$ in \eqref{eq:predictor} by $\lambda_i(x\pm \Delta x)$. We evaluate the source term at $x$, although it is possible to evaluate it at $x \pm \Delta x$ as well. The other steps remain unchanged.

This yields two sets of pairs $(\lambda_i^*, S_i^*)$, one corresponding to centering the predictor step at $x + \Delta x$ and the other at $x - \Delta x$. For each of the sets, $\sum_i \lvert \lambda_i^* \rvert$ is computed, and the one with the largest value is taken. Note that this is slightly different to the original version of the entropy fix, where the decision criterion was applied to every eigenvalue individually. This kind of fix is very simple to implement and has proved itself effective in practice.

\subsection{Near-vacuum treatment} \label{ssec:nearvacuum}

After the entropy fix has selected a set of speeds $\lambda_i^*$, $i = 1, \ldots, M$, it is checked whether any of them is exactly zero. The case \emph{exactly} zero is assumed to only happen if the corresponding initial data at the foot of the characteristic were vacuum. This happens in the vicinity of a dry/wet interface, where at $x$ there is some water height $h > 0$, and non-zero $\lambda_i$-estimates are computed. Then, to compute $\lambda_i^*$, the initial data are evaluated at locations $x + \lambda_i \frac{t}{2}$. If $x$ is close to the dry/wet interface, this location might already be ``on the shore''. Then, one of the $\lambda_i^*$ will be evaluated to zero. In the context of the Euler equations such a setup is a \emph{rarefaction into vacuum}. Obviously, there are no characteristic lines coming from the vacuum: characteristics only can reach out from the wet state. An estimate which locates the footpoint of the characteristic in the dry state is obviously unphysical. In this case, i.e.\ when one of the speeds is zero, we choose to set \emph{all} speeds and values of the source term to zero. 

This ensures that no information flows from vacuum and implies that at this location, the solution remains stationary $q(t, x) = q(0, x)$ over one time step. This means that we prefer to retain the solution at $x$ at this time step and let the rarefaction approach slightly more. In practice, this procedure is found to work well and to avoid spurious oscillations at rarefactions that would otherwise be observed.

Some numerical methods designed to cope with dry/wet interfaces artificially empty cells with little amounts of water. The reason for this is that, on the one hand, little water heights can support spurious hypersonic waves with $v \gg \sqrt{gh}$ which stall the simulation, while on the other hand, most applications of shallow water models lose validity at very low water heights. Physical processes that become dominant in this regime include surface tension, friction, evaporation or percolation. However, setting the water average in a cell to zero is a non-conservative intervention. Here, we suggest to freeze the point values instead, whenever they would be updated to new values smaller than a given threshold $\epsilon$. We choose $\epsilon = 10^{-7}$, because then $\epsilon^2$, the pressure term inside the momentum flux, is at the level of machine precision.  The cell average continues to be updated with frozen fluxes. If after some time the suggested point value lies above the threshold $\epsilon$, it is accepted and used again. As described below, the (frozen) fluxes at the cell interfaces might drain the cell. In this case the adjacent point values are set to $h=0, m=0$.

\subsection{Well-balanced point value update} \label{ssec:wellbalancing}

A static solution to the equations of shallow water \eqref{eq:shallowwaterh}--\eqref{eq:shallowwaterm} is the \emph{lake at rest}:
\begin{align}
 h + b &= \const & m &= 0 \label{eq:lakeatrest}
\end{align}

A numerical method is well-balanced, if it keeps \eqref{eq:lakeatrest} stationary up to machine precision. For Active Flux, achieving a well-balanced evolution consists of two problems: ensuring that the point values remain stationary at the lake-at-rest state, and ensuring the same for the cell averages. The former is discussed now, and the latter -- in Section \ref{ssec:wellbalancingavg}.

It is not surprising that the approximate evolution operator \eqref{eq:evolutionoperator} is not well-balanced. In case the data fulfill \eqref{eq:lakeatrest}, the approximate evolution operator in general will (according to Theorem \ref{thm:evolutionoperator}) yield a time evolution $\mathcal O(t^3)$. However, it can easily become well-balanced with the following strategy inspired by the general well-balancing strategies \cite{pareschi2017,Berberich2019b}. At any location $x_0$ one applies the approximate evolution operator twice: once onto the reconstruction of the actual data, and once onto fictitious initial data of no velocity and a uniformly constant water level $h(x_0) + b(x_0)$. As these fictitious initial data are stationary, any numerical evolution is $\mathcal O(t^3)$. Subtracting it from the evolution of the actual data does not change the order of accuracy of the evolution operator. At the same time, if the data actually do correspond to the lake-at-rest state, then the term subtracted is just the spurious evolution, and the resulting operator keeps the data stationary.

This well-balancing strategy is described now in more detail. Given a fixed location $x_0$, denote the approximate evolution at $x_0$ (including the entropy fix and any other modifications) as $q_1(t, x_0) = (h_1(t, x_0), m_1(t,x_0))$. 

Compute
\begin{align}
 W := h(0, x_0) + b(x_0)
\end{align}
If \eqref{eq:lakeatrest} is true, $W$ is the constant water level.

Apply the approximate evolution operator to initial data $h_0(x) = W - b(x)$, $v_0(x) = 0$ and denote the solution at $x_0$ by $\tilde h(t, x_0)$, $\tilde m(t, x_0) $. Clearly, any actual time evolution is entirely spurious. The well-balanced approximate solution $(h^\text{wb}(t, x_0), m^\text{wb}(t, x_0))$ at $x_0$ is obtained by subtracting the spurious evolution
\begin{align}
 m^\text{wb}(t, x_0) := m_1(t, x_0) - h_1(t, x_0) \frac{\tilde m(t, x_0) }{ \tilde h(t, x_0)} \label{eq:wbevom}\\
 h^\text{wb}(t, x_0) := h_1(t, x_0) - \left(\tilde h(t, x_0) - h(0, x_0)\right)  \label{eq:wbevoh}
\end{align}

\begin{theorem}
 The approximate evolution operator \eqref{eq:wbevom}--\eqref{eq:wbevoh} is well-balanced and has an error $\mathcal O(t^3)$.
\end{theorem}
\begin{proof}
 
If the initial data are a lake at rest, then $h_1(t, x_0) = \tilde h(t, x_0)$ and $m_1(t, x_0) = \tilde m(t, x_0)$ and
\begin{align}
 m^\text{wb}(t, x_0) =&\, \tilde m(t, x_0) - \tilde h(t, x_0) \frac{\tilde m(t, x_0) }{ \tilde h(t, x_0)} = 0\\
 h^\text{wb}(t, x_0) :=&\, \tilde h(t, x_0) - \left(\tilde h(t, x_0) - h(0, x_0)\right)  = h(0, x_0)
\end{align}

The order of accuracy is unchanged because, by the order of accuracy of \eqref{eq:evolutionoperator}, $\tilde h(t, x_0) - h(0, x_0) \in \mathcal O(t^3)$ and $\tilde m(t, x_0) \in \mathcal O(t^3)$.
\end{proof}%

The well-balancing procedure is only applied  
\begin{itemize}
\item if the water height $h$ at $x_0$ is positive, because otherwise the well-balancing tries to induce artificial motion inside vacuum, 
\item if $\tilde h(t, x_0) > 0$, because we need to divide by it, and
\item if the current (absolute) value of the Mach/Froude number $\vert v \vert/\sqrt{gh}$ is less than a threshold (which we choose to be $1$) -- moving waters do not need to be well-balanced. This prevents well-balancing modifications of supersonic phenomena, where clearly no such modification is needed, and saves computational resources. {This decision is taken individually for each point value, such that e.g. a setup of a shock wave running into a lake at rest is resolved correctly.}
\end{itemize}
Note that this well-balancing strategy ignores the details of the evolution operator. It thus can be combined with any approximate evolution operator.

\section{Update of averages} \label{sec:averages}

The update of the averages needs to be such that the resulting method is well-balanced. In the previous chapter it has been shown how to make sure that at the lake-at-rest state the point values are stationary. In Section \ref{ssec:wellbalancingavg} this is used to analyze the necessary condition for the cell averages to remain stationary as well. This completes the well-balancing of Active Flux for shallow water. 

The treatment of dry cells is described in Section \ref{ssec:draining}, and the initialization of cell averages in presence of a dry/wet interface is treated in Section \ref{ssec:initiaizationavg}. 

\subsection{Well-balanced source term quadrature} \label{ssec:wellbalancingavg}

There are several possible choices of a natural quadrature of the source term $- g h \del_x b$. This quadrature, however, needs to be chosen such that the cell averages remain stationary at the lake-at-rest state. Consider therefore
\begin{align}
 h(x) + b(x) &= W \in \mathbb R & v(x) &= 0 \qquad \forall x
\end{align}
and assume that the point value update over one time step evolves this setup exactly, i.e.\ leaves it stationary. In particular assume that the point values remain exactly stationary. Then the average of $h$ automatically remains stationary, but the update of the momentum average reads
\begin{align}
 \del_t \bar m + \frac{\frac12 g \left(h_\text R^2 - h_\text L^2\right)}{\Delta x} = - g \langle \text{discretization of }  h \del_x b \rangle
\end{align}
Using the lake at rest, stationarity of the average requires:
\begin{align}
  - g\langle \text{discr. of }  h \del_x b \rangle\Big \rvert_\text{lake-at-rest} &= \frac{\frac12 g \left(\left(W - b(x_\text R)\right)^2 - (W - b(x_\text L))^2\right)}{\Delta x}  \\
  &\hspace{-1cm}= \frac{\frac12 g \left( b(x_\text R)^2 - b(x_\text L)^2 - 2 W(b(x_\text R) - b(x_\text L)) \right)}{\Delta x} \\
  &\hspace{-1cm}= -g \left(W - \frac{ b(x_\text R)+ b(x_\text L) }{2} \right) \cdot \frac{b(x_\text R) - b(x_\text L)}{\Delta x}  \label{eq:quadraturesourcelakeatrest}
\end{align}

Note that formula \eqref{eq:quadraturesourcelakeatrest} follows uniquely. However, it does not tell the actual discretization of the source term, it only shows what any discretization of the source term needs to reduce to at the lake-at-rest state. A general quadrature formula with this property is given next:

\begin{theorem} \label{thm:quadraturesource}
 The quadrature
 \begin{multline}
  \left( \frac1{12} h^{n+1}_{i+\frac12} + \frac13 h^{n+\frac12}_{i+\frac12} \right )\beta_{1,\text R}
 + \left( \frac1{12} h^{n+1}_{i-\frac12} +\frac13 h^{n+\frac12}_{i-\frac12} \right )\beta_{1,\text L}
 \\-\frac14 h^{n}_{i+\frac12} \beta_{2, \text R}
 -\frac14 h^{n}_{i-\frac12} \beta_{2,\text L}
 +\frac23 h^n_\text{recon}(x_i) \beta_0
 \end{multline}
 with
 \begin{align}
 \beta_{1,\text R} &:= \frac{-4 b_i - b_{i-\frac12} + 5b_{i+\frac12}}{3 \Delta x} &
 \beta_{1,\text L} &:= \frac{4 b_i - 5 b_{i-\frac12} + b_{i+\frac12}}{3 \Delta x}\\
 \beta_{2,\text R} &:= \frac{4 b_i -11 b_{i-\frac12} + 7 b_{i+\frac12}}{9 \Delta x} &
 \beta_{2,\text L} &:= \frac{-4 b_i -7 b_{i-\frac12} + 11 b_{i+\frac12}}{9 \Delta x}\\
 \beta_0 &:= \frac{- b_{i-\frac12} + b_{i+\frac12}}{ \Delta x}
 \end{align}
 approximates
 \begin{align}
  \frac{1}{\Delta x} \frac{1}{\Delta t} \int_0^{\Delta t} \dd t \, \int_{x_{i-\frac12}}^{x_{i+\frac12}} \dd x \, h(x) \del_x b(x)
 \end{align}
 with an error $\mathcal O\left(\Delta t^\alpha \Delta x^\beta\right)$, $\alpha + \beta = 3$ and reduces to \eqref{eq:quadraturesourcelakeatrest} in the lake-at-rest case.
\end{theorem}
\begin{proof}
 Note first that $\beta_0$ and $\beta_{i,\text L/\text R}$ are approximations of $\del_x b(x_i)$. The quadrature formula thus is a linear combination of terms approximating $h \del_x b$. It is exact for quadratic $b$ and for $h$ of the form ($h_{ij} \in \mathbb R$)
 \begin{align}
  h(t, x) = h_{00} + h_{01}x + h_{02} x^2 + t(h_{10} + h_{11}x ) + t^2(h_{20} + h_{21}x )
 \end{align}
 This proves the order of approximation. When all values of $h$ are stationary, then the formula reduces to
  \begin{align}
  &  h^{n}_{i+\frac12} \left( \frac5{12} \beta_{1,\text R}  -\frac14  \beta_{2, \text R} \right )
    + h^{n}_{i-\frac12} \left( \frac5{12}  \beta_{1,\text L}  -\frac14  \beta_{2,\text L} \right )
 +\frac23 h^n_\text{recon}(x_i) \beta_0\\
 &= h^{n}_{i+\frac12}   \frac{-24 b_i + 6 b_{i-\frac12} + 18 b_{i+\frac12}}{36 \Delta x} 
    + h^{n}_{i-\frac12} \frac{24 b_i - 18 b_{i-\frac12} - 6 b_{i+\frac12} }{36 \Delta x} 
 +\frac23 h^n_\text{recon}(x_i) \beta_0\nonumber
 \end{align}
 which gives -- when the lake is at rest --
  \begin{align}
     (W - b_{i+\frac12})  \frac{-24 b_i + 6 b_{i-\frac12} + 18 b_{i+\frac12}}{36 \Delta x} &\\
    + (W - b_{i-\frac12})  \frac{24 b_i - 18 b_{i-\frac12} - 6 b_{i+\frac12} }{36 \Delta x} &\nonumber
 \\+\frac23 (W - b_i) \frac{- b_{i-\frac12} + b_{i+\frac12}}{ \Delta x}
 &= \nonumber
  W   \frac{b_{i+\frac12} -  b_{i-\frac12}  }{\Delta x} 
    -   \frac{  b_{i+\frac12}^2 - b_{i-\frac12}^2  }{2 \Delta x} 
 \end{align}
\end{proof}

\emph{Note}: If $h$ is spatially constant at all times, the quadrature formula reduces to Simpson's rule in time.

In order to maintain the lake at rest exactly, the quadrature needs to give the exact value at partially wet cells. At partially flooded cells, the location of the shore is computed as $x^*$ during the reconstruction procedure. If the dry state is on the left (right) of the cell, instead of $h^n_\text{recon}(x_i)$ we therefore suggest to use $h^n_\text{recon}\left(\frac{x^* + \frac{\Delta x}{2}}{2}\right)$ or $h^n_\text{recon}\left(\frac{-\frac{\Delta x}{2} + x^*}{2}\right)$.

\subsection{Draining} \label{ssec:draining}

Before updating the cell averages there is a trial run of computing
\begin{align}
 \bar h_i - \Delta t \frac{\hat f^h_{i+\frac12} - \hat f^h_{i-\frac12}}{\Delta x}
\end{align}
for all cells, where $\hat f_{i+\frac12}^h$ denotes the numerical flux of $h$, obtained as quadrature \eqref{eq:fluxquadrature}. (Recall that there is no source term in the equation for $h$.) During this run, it is checked whether any cell would acquire a negative\footnote{We only check whether the water height would be less than $-10^{-14}$, because cells with a water average $\bar h\in[-10^{-14},10^{-14}]$ are emptied completely to avoid dealing with machine errors.} average water height. In this case, following \cite{bollermann11}, we compute the draining time $\Delta t_\text{drain}$ via
\begin{align}
 \bar h_i - \Delta t_\text{drain} \frac{\hat f^h_{i+\frac12} - \hat f^h_{i-\frac12}}{\Delta x} = 0
\end{align}
and reduce all the fluxes across the edges of this cell as follows:
\begin{align}
 \hat f_{i+\frac12} \mapsto \frac{\Delta t_\text{drain}}{\Delta t} \hat f_{i+\frac12}
\end{align}
Note that this applies to all fluxes, also the momentum flux. This procedure makes sure that maintaining non-negative water height averages does not reduce the CFL condition.

After finishing the loop, the loop is reentered again, until no cell requires a change. This is because a modification of the flux might entail a change of the flux balance in a neighbouring cell. An efficient way of implementing this is to concentrate on the neighbours of a cell which has been identified to be drained. Note however, that the recursion is sure to stop, as the fluxes are constantly decreasing in absolute value, and the current averages are all non-negative.

Having corrected the fluxes, first, the flux difference is used to update the cell averages. This way, cells that would have acquired negative averages of water height $h$ at most become empty. For all cells, for which the fluxes have led to an empty cell, the adjacent point values are set to zero water height and zero momentum. This ensures that a continuous reconstruction remains possible. Such a procedure also makes sense intuitively.

Finally, in all cells with non-zero water height averages, the source term quadrature according to Theorem \ref{thm:quadraturesource} (multiplied by $\Delta t$) is added to the momentum average.

Note that in our method, cells are only considered dry, when the average water height is machine zero (i.e. $10^{-14}$). This is in contrast to other works where a cell is considered dry if the water height is smaller than some larger threshold.

\subsection{Initialization of cell averages} \label{ssec:initiaizationavg}

The initialization of cell averages for Active Flux requires a point value e.g.\ at the cell center. Together with the point values at the cell interfaces, the cell average is initialized via Simpson's rule. Inside a partially flooded cell this gives a wrong result, usually not dramatic, but spoiling well-balancing. In order to initialize the cell average correctly in this case, we use a ``shore search''. We assume that the bottom topography is sufficiently resolved, such that any computational cell initially contains at most one dry/wet interface. Its location is found via a binary search to machine precision. If $h_\text R = 0$, and the shore is found to be located at $x = x_\text s$, then the average of both the water height and momentum is computed as
\begin{align}
 \bar q = \frac16 \left( q_\text L + 4 q_0\left( \frac{x_\text s + x_\text L}{2}  \right) + q_0(x_\text s)\right) \frac{x_\text S - x_\text L}{\Delta x}
\end{align}
where $x_\text L$ is the location of the left interface of the cell, and $q_0$ are the initial data. An analogous calculation applies in the mirrored situation. Note that this procedure occurs only in the initialization and thus is not time-critical.

\section{Numerical examples} \label{sec:numerical}

In all the tests we use $g = 9.812$. 

\subsection{Stationary solutions and well-balancing}

\begin{figure}
 \centering
 \includegraphics[width=\textwidth]{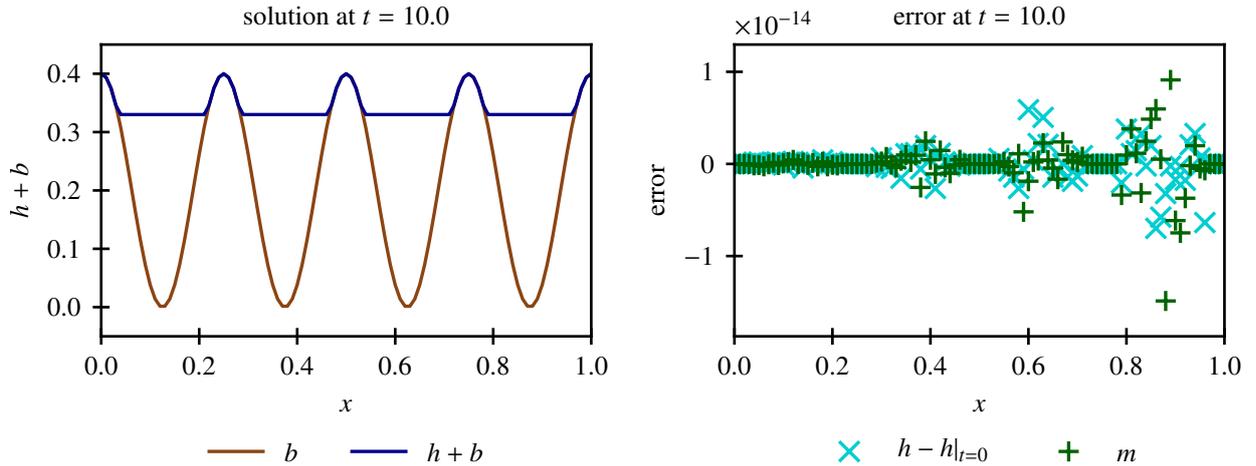}
 \caption{Demonstration of the well-balanced property in presence of partially dry cells. \emph{Left}: Setup with four lakes at rest. Point values of $h+b$ are shown. \emph{Right}: Errors of the point values of the numerical solution at $t=10$.}
 \label{fig:wellbalancing}
\end{figure}

The first test assesses the ability of the numerical method to preserve a lake at rest to machine precision even when partially dry cells are involved. On $[0,1]$ discretized with $100$ cells, we use
\begin{align}
 b(x) = 0.2 (1 + \cos(8 \pi x))
\end{align}
and $h(0, x) = \max(0, 0.33 - b(x))$, $m(0, x) = 0$ (see Figure \ref{fig:wellbalancing}). We find that this setup is kept stationary to machine precision.

{%
Figure \ref{fig:wellbalancingperturbation} shows the evolution of a Gaussian water height perturbation in one of the lakes.

\begin{figure} 
 \centering
 \includegraphics[width=\textwidth]{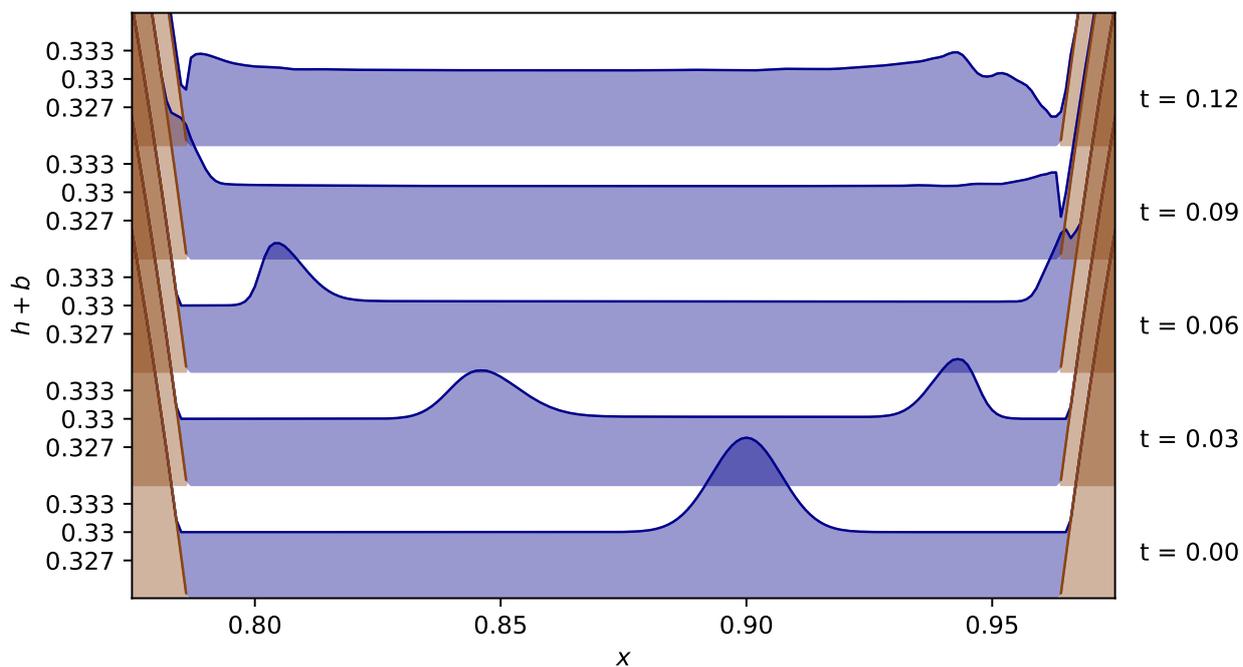}
 \caption{Evolution of a Gaussian pertubation in one of the lakes of Figure \ref{fig:wellbalancing}. Point values of $h+b$ are shown; $\Delta x = 10^{-3}$.}
 \label{fig:wellbalancingperturbation}
\end{figure}
}

\subsection{Convergence}

\begin{figure} 
 \centering
 \includegraphics[width=\textwidth]{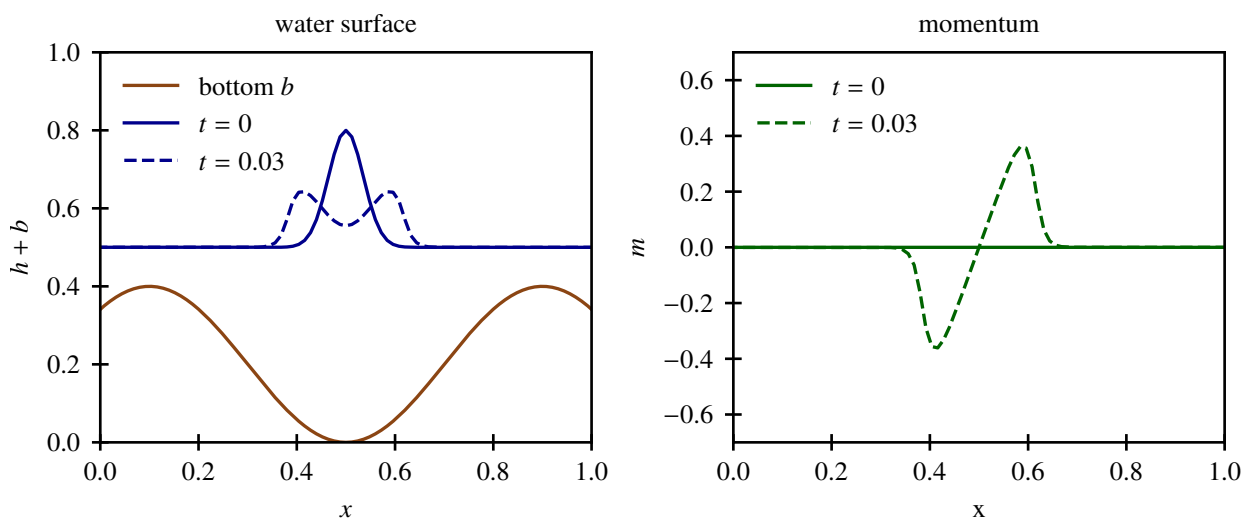}
 \caption{Setup for the convergence test and its numerical evolution on a grid with 200 cells. \emph{Left}: Water level $h+b$. \emph{Right}: Momentum.}
 \label{fig:convergencesetup}
\end{figure}

\begin{figure}
 \centering
 \includegraphics[width=0.7\textwidth]{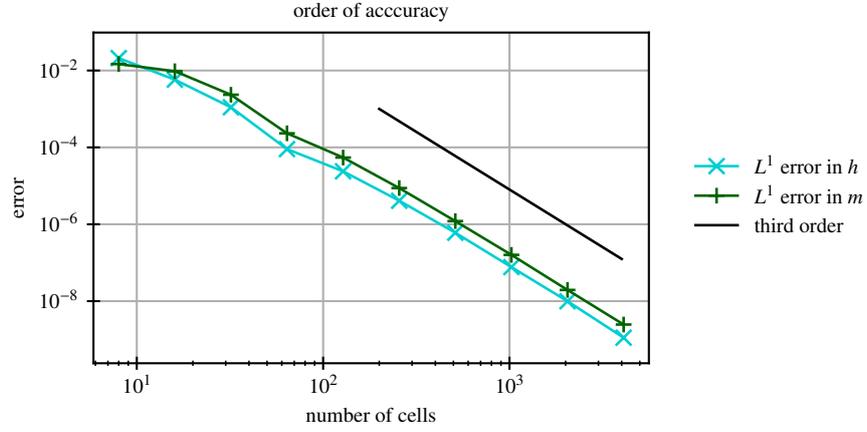}
 \caption{Convergence of a Gaussian wave on cosine-shaped bottom. The $L^1$ error of the point values is shown.}
 \label{fig:convergence}
\end{figure}

To demonstrate third-order accuracy of the numerical method we use
\begin{align}
 b(x) &= 0.2 (1 + \cos(6 \pi x))\\
 h(0, x) + b(x) &= 0.5 + 0.3 \exp\left(- \frac{(x-0.5)^2}{0.05^2}\right)
\end{align}
and $m(0,x) = 0$. This setup is evolved until $t=0.03$ with a CFL of $0.7$ (see Figure \ref{fig:convergencesetup}). The reconstruction is parabolic without any limiting or positivity-preserving modification. However, well-balancing remains active. The $L^1$ error (shown in Figure \ref{fig:convergence} for various resolutions) is computed using the solution obtained on a grid with $16384 = 2^{14}$ cells as reference. The domain is $[0, 1]$ with periodic boundaries.

Figure \ref{fig:convergence} {and Table \ref{tab:convergence} demonstrate} that the numerical method attains third order of accuracy, as expected.

\begin{table}
 \centering
 \begin{tabular}{c||c|c||c|c}
grid size & $L^1$ error $h$ & order & $L^1$ error $m$ & order\\ \hline
64	&$2.32242 \cdot 10^{-4}$&		&$9.02165 \cdot 10^{-5}$&	\\
128	&$5.44272 \cdot 10^{-5}$&	2.09	&$2.38423 \cdot 10^{-5}$&	1.92\\
256	&$8.73774 \cdot 10^{-6}$&	2.64	&$4.09151 \cdot 10^{-6}$&	2.54\\
512	&$1.20775 \cdot 10^{-6}$&	2.85	&$6.01095 \cdot 10^{-7}$&	2.77\\
1024	&$1.59633 \cdot 10^{-7}$&	2.92	&$7.66888 \cdot 10^{-8}$&	2.97\\
2048	&$1.95048 \cdot 10^{-8}$&	3.03	&$9.92063 \cdot 10^{-9}$&	2.95\\
4096	&$2.45793 \cdot 10^{-9}$&	2.99	&$1.11746 \cdot 10^{-9}$&	3.15\\
 \end{tabular}
 \caption{$L^1$ errors and numerical order of accuracy for the convergence test of Figure \ref{fig:convergencesetup}.}
 \label{tab:convergence}
\end{table}

{%

\subsection{Accuracy test}

This accuracy test has been suggested in \cite{bouchut04} (Chapter 6) and involves a continuous solution, with discontinuous derivatives. The bottom topography is
\begin{align}
 b(x) = \max\left(0,  0.48 \left( 1 - \left(  \frac{x-20}{4} \right)^2 \right )  \right)
\end{align}
and the initial data are $h(0,x) = 4$, $m(0,x) = 10$. The numerical results at $t=1$ with periodic boundaries are shown in Figure \ref{fig:accuracybouchut}. Table \ref{tab:accuracybouchut} summarizes the $L^1$ errors upon mesh refinement, computed with respect to a highly resolved solution on a grid of 2048 cells. The presence of discontinuities in the derivatives reduces the order of accuracy; for comparison, second-order methods in \cite{bouchut04} attain an experimental order of accuracy between 1.25 and 1.5.

\begin{figure} 
 \centering
 \includegraphics[width=\textwidth]{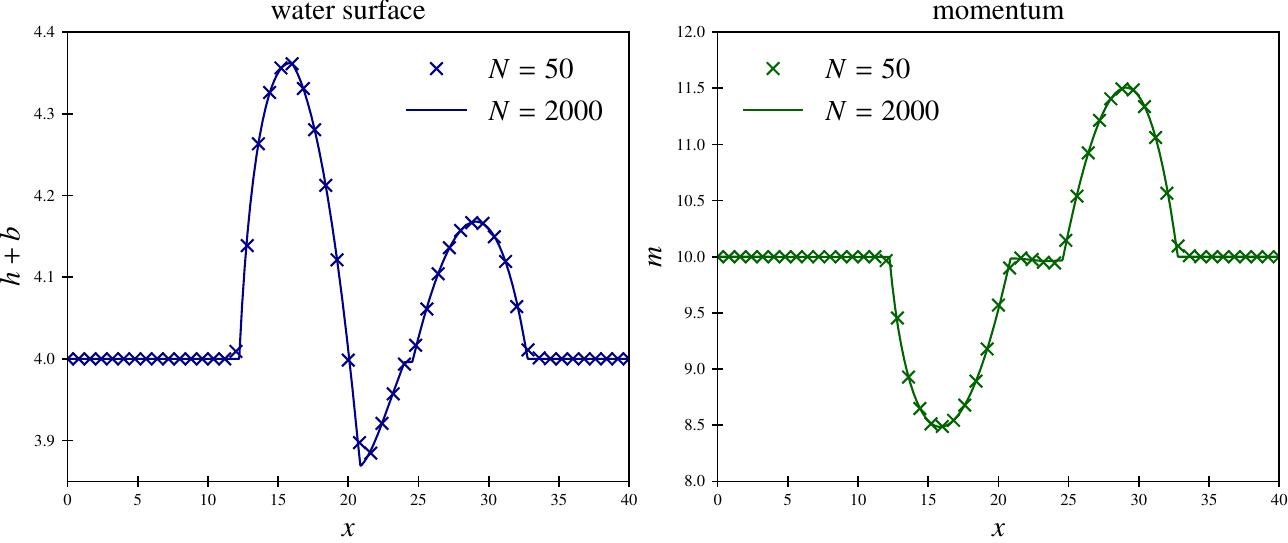}
 \caption{Accuracy test from \cite{bouchut04}, showing results (point values) on a grid of 50 cells (crosses) and 2000 cells (solid line). The CFL number is 0.7. \emph{Left}: Water level $h+b$. \emph{Right}: Momentum $m$.}
 \label{fig:accuracybouchut}
\end{figure}

\begin{table} 
 \centering
 \begin{tabular}{c||c|c||c|c}
 grid size & $L^1$ error $m$ & order & $L^1$ error $h$ & order\\\hline
  8&4.36965206	& 	&0.73838181	&\\
16&	1.51920049	&1.52	&0.29898920&	1.30\\
32&	0.39594621	&1.94	&0.07952090&	1.91\\
64&	0.16409114	&1.27	&0.02966787&	1.42\\
128&	0.08921846	&0.88&	0.01707322&	0.80\\
256	&0.02911097	&1.62	&0.00568257	&1.59\\
512	&0.00898213	&1.70	&0.00175999	&1.69\\
1024&	0.00258718	&1.80&	0.00052705&	1.74
 \end{tabular}
 \caption{$L^1$ errors and numerical order of accuracy for the test of Figure \ref{fig:accuracybouchut}.}
 \label{tab:accuracybouchut}
\end{table}

}

\subsection{Riemann problems}

\begin{figure}
 \centering
 \includegraphics[width=\textwidth]{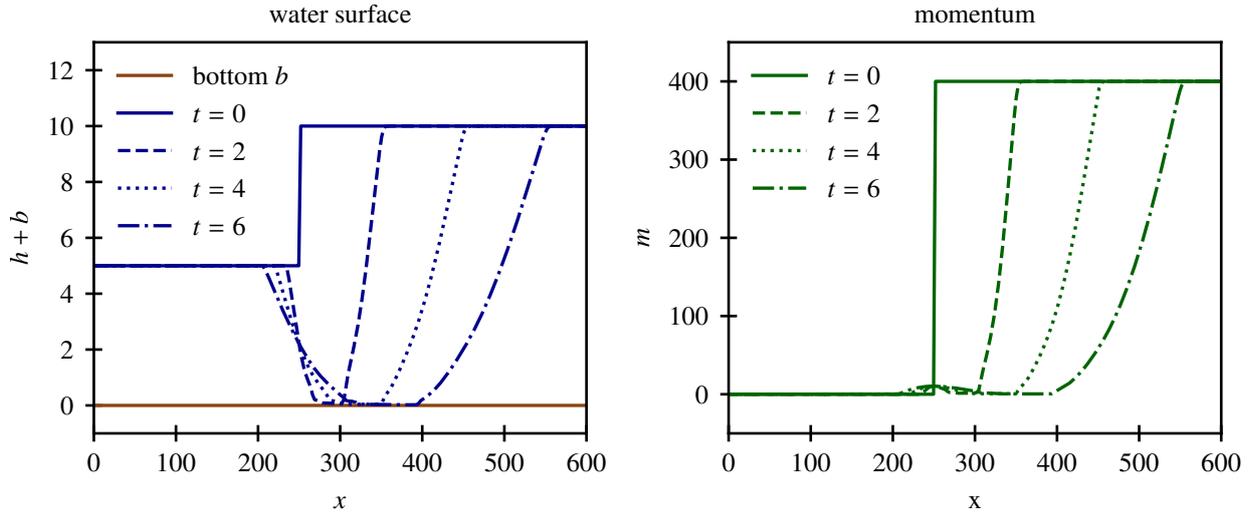}
 \caption{Double-rarefaction Riemann problem from \cite{xing11}. Point values are shown at different times. \emph{Left}: Water height. \emph{Right}: Momentum.}
 \label{fig:doublerarefactionriemann}
\end{figure}

The first Riemann problem (from \cite{xing11}) involves a strong rarefaction that produces vacuum. The bottom topography is flat in this case. This is a good test for both the ability of the code to deal with discontinuous setups and for its non-negativity preserving properties. We solve it on a grid of 250 cells covering $[0,600]$ with a CFL number of 0.7. As can be seen in Figure \ref{fig:doublerarefactionriemann}, the Active Flux method does not have difficulty resolving the solution, see \cite{xing11} for comparison with results obtained using a WENO method.

The bottom topography is projected onto globally continuous piecewise polynomials. Approximating a discontinuous bottom with steep gradients allows to compute Riemann problems that can be compared to those of discontinuous bottom studied in \cite{chinnayya04}. Next, we thus use the following setup:
\begin{align}
 b(x) &= \begin{cases} B & x < 0.5 \\
		B\left(1  - \frac{x - 0.5}{\epsilon/2}\right) & x < 0.5 + \epsilon/2 \\
		0 & \text{else} \end{cases}\\
 h(0,x) &= \begin{cases} 3 & x < 0.5 \\ 4 & \text{else} \end{cases}
\end{align}
and $m(0,x) = 0$ with $\epsilon = 0.01$ and $B \in \{2, 4, 50 \}$, corresponding to Figures 10, 12 and 15 in \cite{chinnayya04}. To compare our results with those of \cite{chinnayya04}, we evolve the setup until $t=0.048$. A CFL number of $0.7$ is used, and $\Delta x = 10^{-3}$ for $B \in \{2, 4\}$ and $\Delta x = 0.8 \cdot 10^{-3}$ for $B = 50$.

\begin{figure} 
 \centering
 \includegraphics[width=\textwidth]{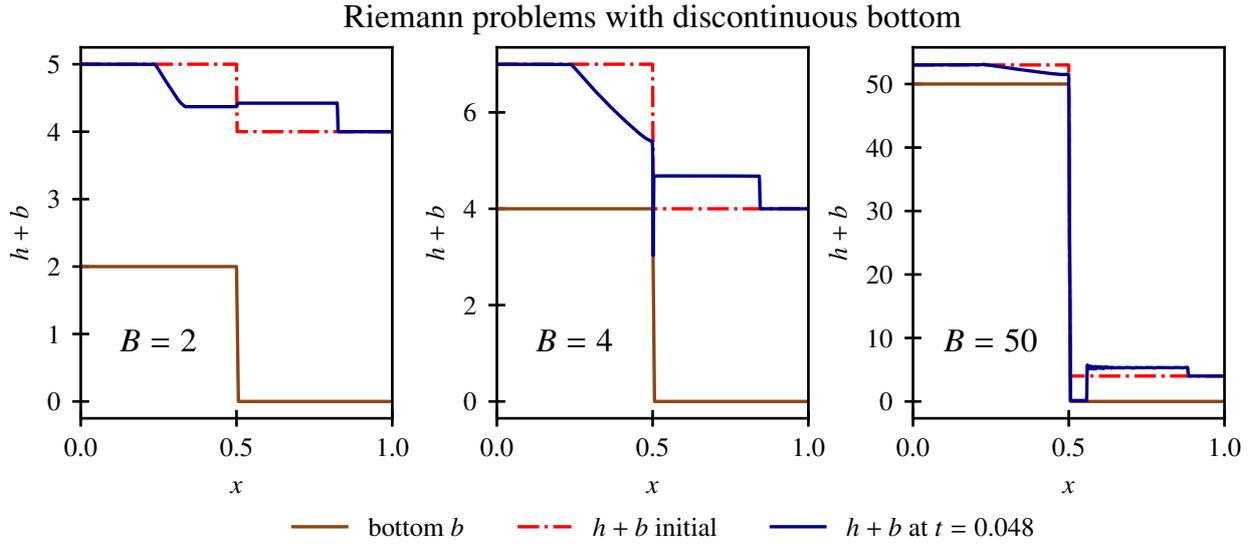}
 \caption{Test problems from \cite{chinnayya04}, but solved using a regularized (continuous) bottom. Point values of $h+b$ are shown at $t=0.048$.}
 \label{fig:seguin}
\end{figure}

\begin{figure}
 \centering
 \includegraphics[width=\textwidth]{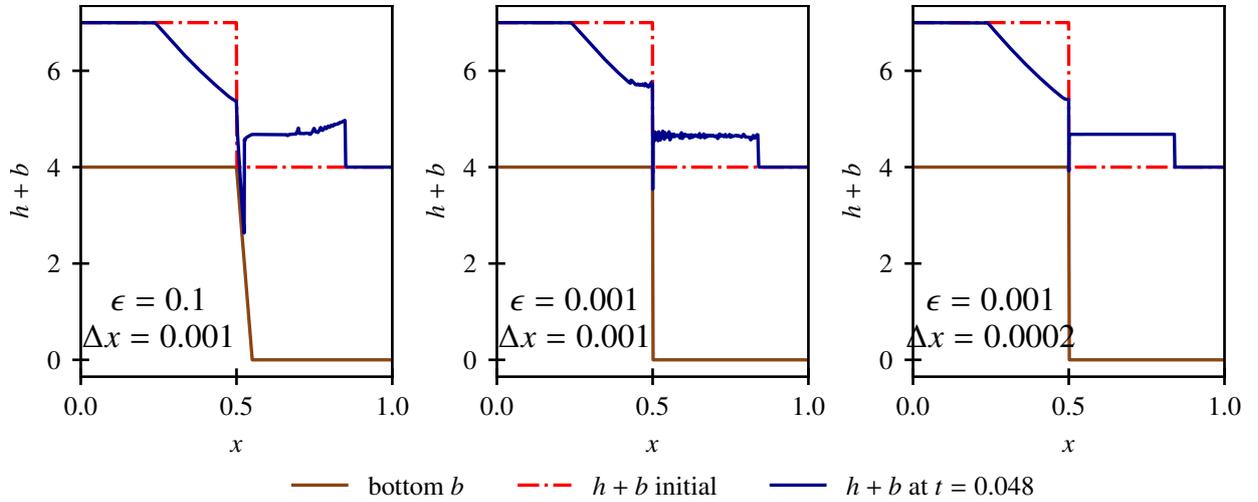}
 \caption{Study of the influence of $\epsilon$ for test problem with $B=4$. Additionally to the case $\epsilon = 0.01$ with $\Delta x = 10^{-3}$ of Figure \ref{fig:seguin} (center) the following cases are considered. \emph{Left}: $\epsilon = 0.1$ on a grid with $\Delta x = 10^{-3}$. One observes that the solution carries a visible imprint of the jump regularization by a steep slope. The oscillations disappear upon grid refinement. \emph{Center}: $\epsilon = 0.001$ with $\Delta x = 10^{-3}$. One observes the rarefaction to be resolved badly. \emph{Right}: $\epsilon = 0.001$ with $\Delta x = 2 \cdot 10^{-4}$. The correct solution is obtained again.}
 \label{fig:seguinepsilon}
\end{figure}

The numerical results are presented in Figure \ref{fig:seguin}. These challenging setups are solved without difficulty. Observe that the case $B=4$ is transonic with two stationary waves (1-shock and 0-wave) superimposed at the location of the bottom discontinuity (see \cite{chinnayya04} for more details on this so-called \emph{resonant} case, and in particular Figure 14 therein). With the regularized bottom, these two waves can be distinguished, and the spike that we observe is not a numerical artefact. It is present because of the finite regularization $\epsilon > 0$. {Figure \ref{fig:seguinepsilon} demonstrates how upon choosing a smaller $\epsilon$ one needs to also reduce $\Delta x$. The two waves are then a few computational cells apart from each other.} 

\subsection{Parabolic bowl}

\begin{figure}
	\centering
	\includegraphics[width=\textwidth]{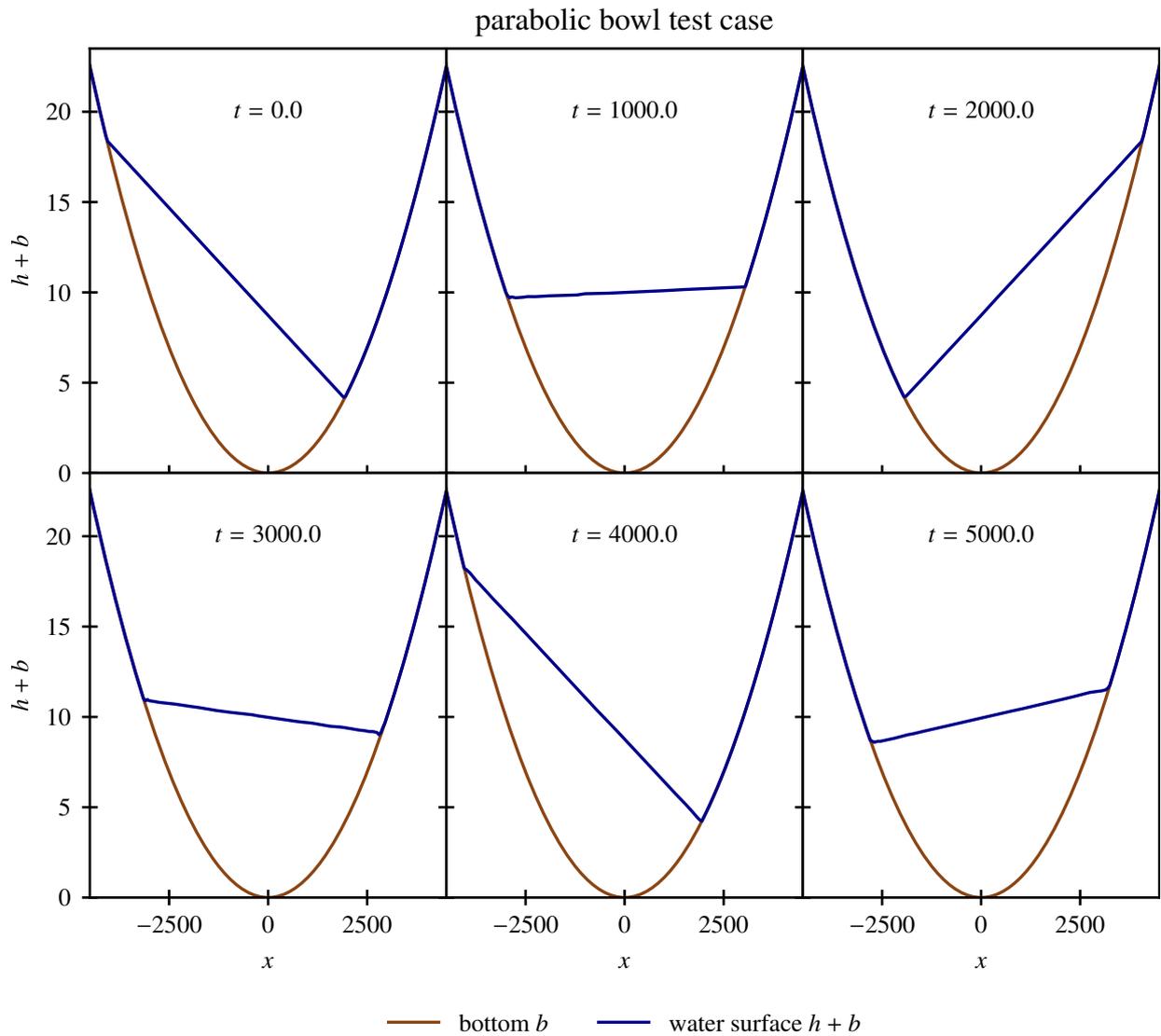}
	\caption{Sloshing water in a parabolic bowl. Point values of $h+b$ are shown at $t = 0, 1000, \ldots, 5000$.}
	\label{fig:parabolicbowl}
\end{figure}

The shallow water equations admit a periodic-in-time solution (\cite{thacker81}) with a parabolic bottom $b = \left(\frac{x}{x_0}\right)^2$ and a water level linear in space:
\begin{align}
 h(t, x) + b(x) &= H_0 - \frac{v_\text{max}^2 }{4g} - \frac{v_\text{max}^2}{4 g} \cos(2\Omega t) - \sqrt\frac{2}{g} \frac{v_\text{max}}{x_0} \cos(\Omega t) x \\
 v(t,x) &=  v_\text{max} \sin(\Omega t)
\end{align}
where $x_0, H_0, v_\text{max}$ are arbitrary constants and
\begin{align}
 \Omega^2 = \frac{2 g}{x_0^2}
\end{align}
The shoreline is given by
\begin{align}
 x_\text{shore}(t) &= \frac{x_0}{2}\left(\pm 2 \sqrt{H_0} - \sqrt\frac{2}{g} v_0 \cos(\Omega x) \right) \label{eq:shorelineevo}
\end{align}
Following \cite{xing10,vater14} we choose $v_\text{max} = 5$, $x_0 = 300\sqrt{10}$ and $H_0 = 10$. The grid consists of 200 cells covering $[-5000,5000]$, and the simulation runs until $t=5 \cdot 10^3$ with CFL $=0.7$.
$H_0$ is the water level when the surface is horizontal. During the oscillation, the dry/wet boundary reaches a height of
\begin{align}
 \frac{\left( \sqrt{2} v_\text{max} + 2 \sqrt{g H_0} \right )^2}{4g} \simeq 18.41
\end{align}
The oscillation period is $\frac{2\pi}{\Omega} = 1345.57$. Figure \ref{fig:parabolicbowl} shows $h+b$ at times $t=0$, $1000$, $\ldots$, $5000$. One observes that the numerical method is able to maintain the linear water surface. Moreover, Figure \ref{fig:parabolicbowlshore} shows the time evolution of the shore, which matches closely the exact evolution \eqref{eq:shorelineevo}. {Figure \ref{fig:parabolicbowlerror} shows that the order of accuracy for this setup is about 2, which probably is due to the presence of the shore, and thus reduced regularity of the solution.}
%

\begin{figure}
 \centering
 \includegraphics[width=\textwidth]{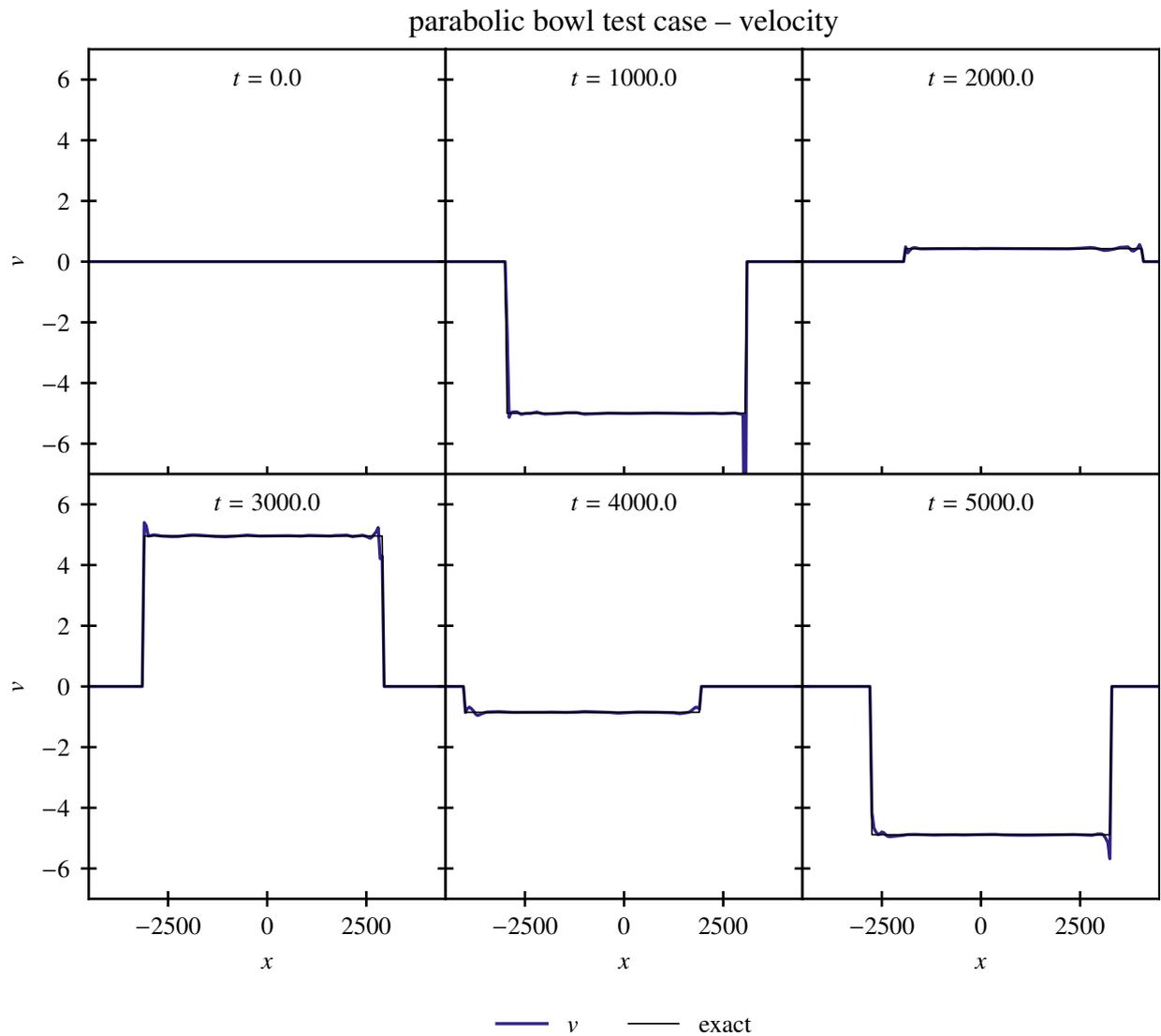} 
 \caption{{The velocity of sloshing water in a parabolic bowl is shown together with the exact solution.}}
 \label{fig:parabolicbowlshorevelocity}
\end{figure}

\begin{figure}
 \centering
 \includegraphics[width=\textwidth]{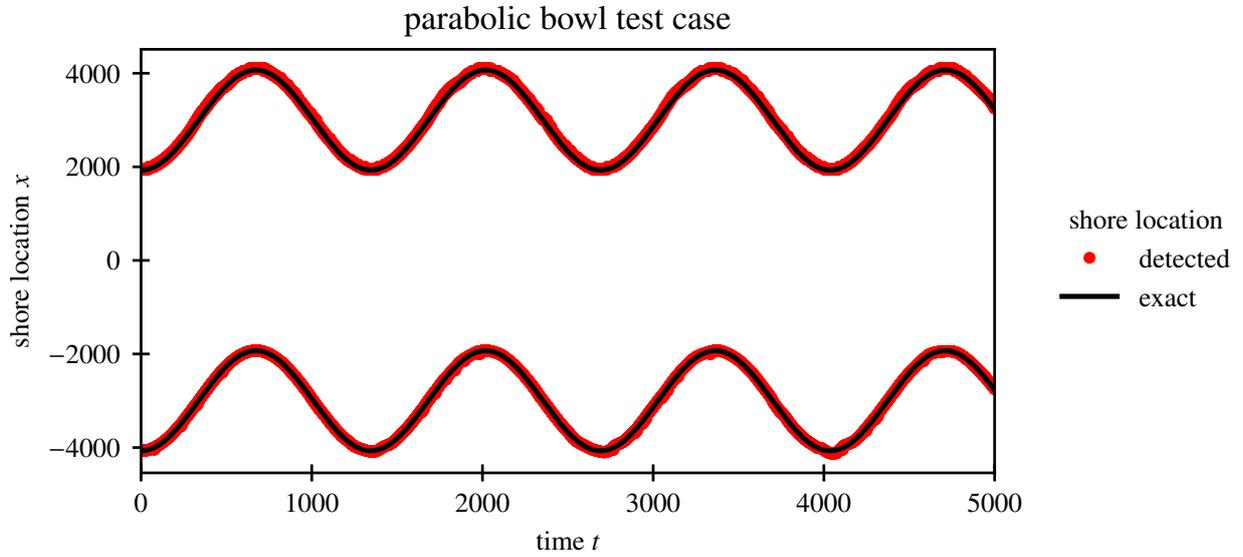} 
 \caption{The location of the shore for sloshing water in a parabolic bowl (points) is shown together with the exact solution (solid line) as a function of time.}
 \label{fig:parabolicbowlshore}
\end{figure}

\begin{figure} 
 \centering
 \includegraphics[width=0.7\textwidth]{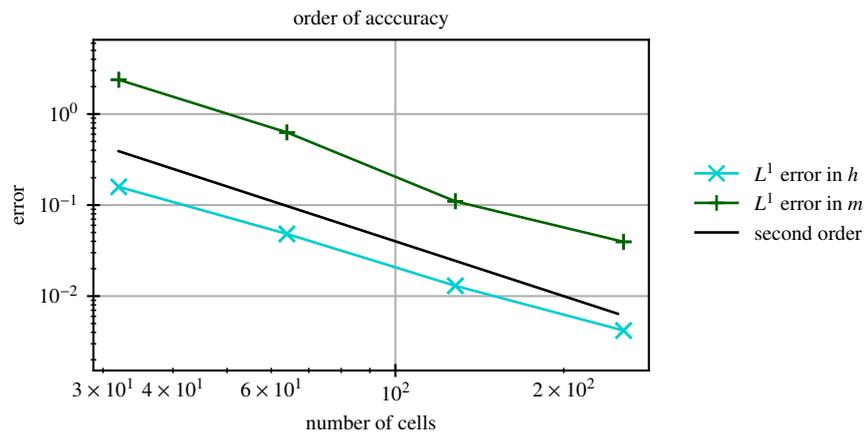}
 \caption{{Error of the numerical solution (point values) computed at $t=5000$ for different grid sizes.}}
 \label{fig:parabolicbowlerror}
\end{figure}

\subsection{Double rarefaction}

This test is from \cite{gallouet03,bouchut04} (Section 6, Test 3). On a domain $[0,25]$, the topography and initial data are given as follows:

\begin{minipage}[l]{0.4\textwidth}
\begin{align*}
 b(x) &= \begin{cases} 
 0 & x < \frac{25}{3} - \frac{\epsilon}{2} \\
 1 & x \in \left[ \frac{25}{3},\frac{25}{2} \right] \\ 
 0 & x > \frac{25}{2} + \frac{\epsilon}{2}\\
 \text{linear} & \text{else} 
 \end{cases} \end{align*}
 \end{minipage}\hfill
 \begin{minipage}[r]{0.55\textwidth} \begin{align*}
 h(0,x) &= 10 - b(x) \\
 m(0,x) &= \begin{cases} -350 & x < \frac{50}{3}\\ 350 & \text{else} \end{cases}
\end{align*}
\end{minipage}

As our method employs a continuous bottom, a comparison with a setup employing a discontinuous bottom is possible by approximating the discontinuity by a steep gradient. The non-constant parts of the bottom are such that $b$ is continuous; their support has size $\epsilon$.

\begin{figure}  
 \includegraphics[width=\textwidth]{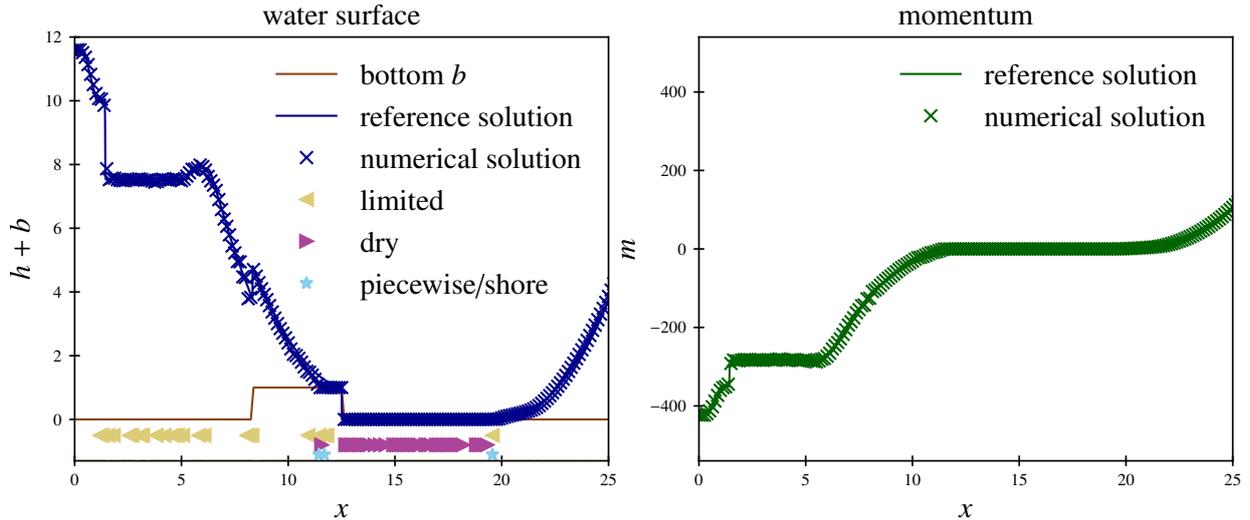}
 \caption{Double rarefaction test from \cite{gallouet03,bouchut04}. The computation is performed on a larger grid to avoid the influence of boundaries; $\Delta x = \frac{25}{200}$ and the CFL number is $0.7$. The point values of the water level $h + b$ (\emph{left}) and of the momentum $m$ (\emph{right}) are shown at time $t=0.25$ for $\epsilon=0.01$, together with a reference solution (solid line) computed on a grid with $\Delta x = \frac{25}{2000}$. In the left plot, symbols indicate the kind of reconstruction in the corresponding cell (limited, dry, or piecewise defined because half-wet). The vertical location of the symbols has no importance.}
 \label{fig:doublerarefaction}
\end{figure}

Numerical results are shown in Figure \ref{fig:doublerarefaction}, and the reader is referred to \cite{bouchut04} for comparison with other schemes. In the figure, the reconstruction of the various cells is also shown, indicating in which cells the limiter is active, and which cells are reconstructed as dry.

\subsection{Transcritical flow with shock}

\begin{figure}  
 \centering
 \includegraphics[width=0.7\textwidth]{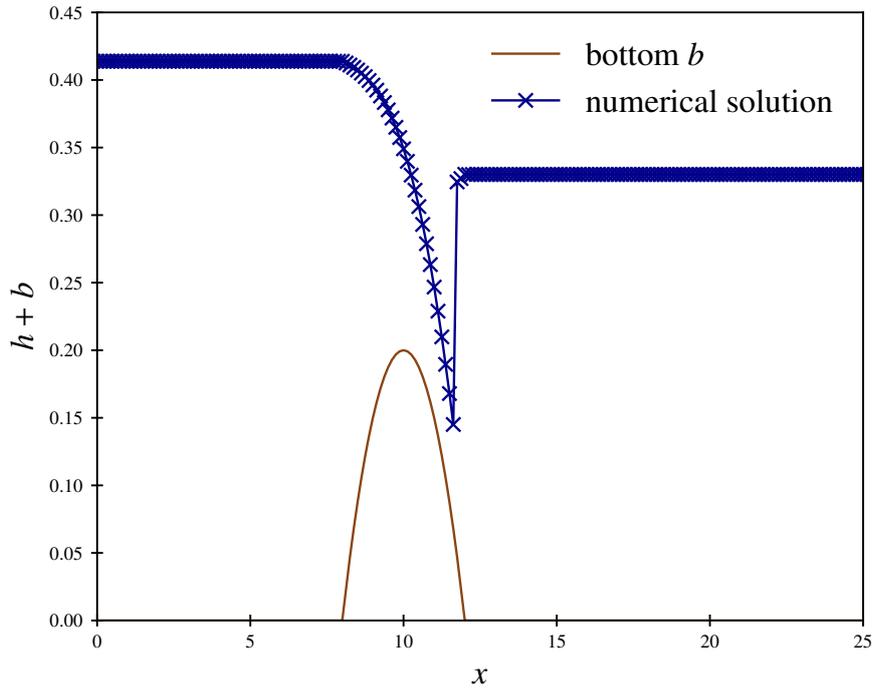}
 \caption{Transcritical shock test; $\Delta x = \frac{25}{200}$ and the CFL number is $0.7$. The point values of the water level $h + b$ are shown at time $t=50$.}
 \label{fig:transcritical}
\end{figure}

For this test, on a domain covering $[0,25]$ the bottom topography is chosen as
\begin{align}
 b(x) = \max(0, 0.2 - 0.05 (x-10)^2)
\end{align}
with initial data $h(0,x) = 0.33$, $m(0,x) = 0.18$ and Dirichlet boundary conditions (see e.g. \cite{bouchut04,gallouet03}). The numerical simulation stationarizes on a transcritical shock, which can be considered a model for hydraulic jumps readily observable in rivers behind underwater obstacles. The results are shown in Figure \ref{fig:transcritical}; one fails to observe any nonentropic features or oscillations.

\subsection{Wave run-up on a sloped beach}

As a last test case we compute a wave run-up onto a plane beach. Following \cite{vater14}, we use the setup of benchmark problem 1 from the Third International Workshop on Long-Wave Runup Models\footnote{June 17-18 2004, Wrigley Marine Science Center, Catalina Island, California, \url{http://isec.nacse.org/workshop/2004_cornell/}}. This benchmark is provided with both the initial data and the analytical solutions (using \cite{carrier03}) at $t\in \{160,175,220\}$ and the time evolution of the shore.

\begin{figure}
 \centering
 \includegraphics[width=\textwidth]{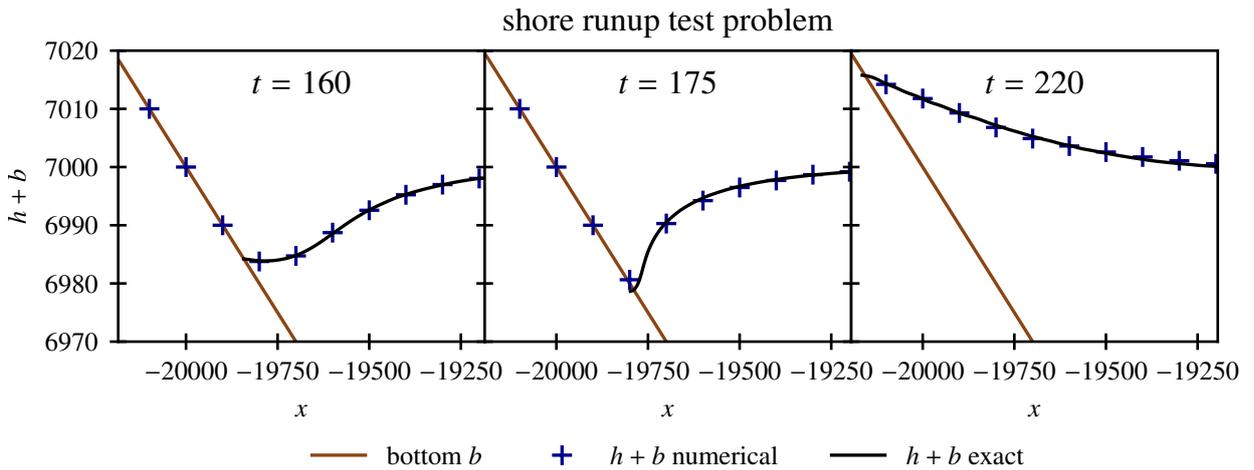}
 \caption{Wave run-up benchmark problem. Point values of $h+b$ at $t = 160, 175, 220$ are shown together with the analytical solution (solid line).}
 \label{fig:tsunami}
\end{figure}

\begin{figure}
 \centering
 \includegraphics[width=\textwidth]{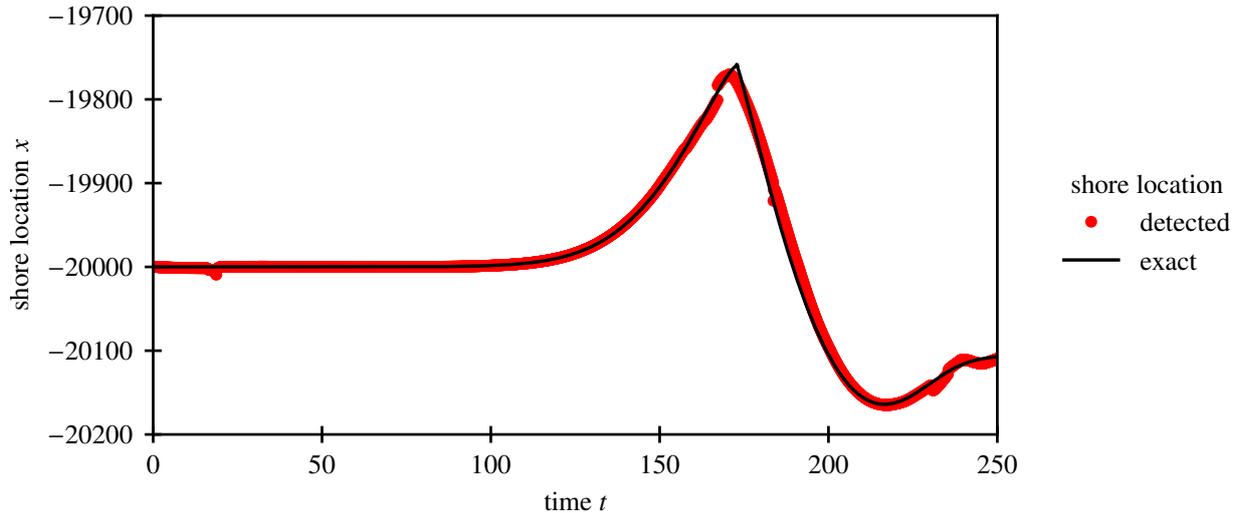}
 \caption{Tsunami run-up benchmark problem. The time evolution of the shore (crosses) is shown together with the analytical solution (solid line). The oscillations visible in the receding wave disappear upon grid refinement.}
 \label{fig:tsunamishore}
\end{figure}

Our results are shown in Figure \ref{fig:tsunami}. We use $500$ cells covering $[-25000,25000]$ and a CFL of 0.7. Figure \ref{fig:tsunamishore} shows the temporal evolution of the shore. 

\section{Conclusions and outlook}

Active Flux is a novel numerical method for systems of hyperbolic conservation laws and is an extension of the finite volume method. Additionally to the cell average, it actively evolves point values at cell interfaces. This yields a high order numerical method with a continuous reconstruction -- an approach very different from, e.g., Discontinuous Galerkin methods. Instead of a short-time evolution of discontinuous initial data (Riemann solver), a short-time evolution of continuous data is used, which equally provides upwinding for stability. Instead of directly updating the average, the point values are updated first, and they are then used to compute the numerical flux via quadrature. Despite a continuous reconstruction, Active Flux is able to accurately resolve Riemann Problems.

The continuous reconstruction is believed to be one of the reasons for the naturally arising structure preservation properties of Active Flux. Whereas even exact Riemann solvers often do not yield a structure preserving method, Active Flux has so far been observed to be structure preserving upon usage of an exact evolution operator. Of course, in practice, approximate evolution operators are needed, but the necessary modification to achieve structure preservation is small, and so far seems to be easy to find. These structure preservation properties make Active Flux a promising candidate for new, accurate, high-order numerical methods.

As a stepping stone towards understanding and demonstrating the abilities of Active Flux, this paper applies it to the equations of shallow water, and achieves a well-balanced third order method that is able to deal with dry states. This brings Active Flux closer to complex applications beyond pure academic interest.

Because Active Flux does not use continuous reconstructions, its implementation for the shallow water equations as presented in this paper also uses a continuous bottom topography. Piecewise parabolic bottom topography, together with a parabolic reconstruction of the water height in each cell, allows to reconstruct the lake at rest exactly. The continuous bottom topography also does not necessitate to deal with products of discontinuous functions in the source term. However, the treatment of the dry/wet interface becomes more complicated.

This paper shows an example of how Active Flux can be successfully applied to the one-dimensional shallow water equations and account for wetting and drying. After having been suggested as early as in \cite{vanleer77}, Active Flux has been developed further in \cite{eymann13} because of its structure preserving properties in particular in multiple spatial dimensions. As this paper shows, it is an original and competitive numerical method for one-dimensional problems already. But clearly, many applications require the availability of a multi-dimensional numerical method. From the beginning, Active Flux was designed to be particularly accurate when applied to multiple spatial dimensions. To extend it in, say, a directionally split fashion would be against the philosophy of Active Flux and might possibly spoil structure preservation. Indeed, past efforts have focused on the observation that characteristics need to be conceptually replaced by characteristic cones in multiple spatial dimensions, and that reflecting this fact in the numerical method is crucial for structure preservation. Whereas this has been studied extensively for linear multi-dimensional hyperbolic systems (e.g. in \cite{barsukow18activeflux}), multi-dimensional evolution operators for nonlinear systems remain subject of future work. However, it is clear that most of the strategies developed in this paper (e.g. those concerning well-balancing and the treatment of dry states) would immediately carry over to the multi-dimensional method. Future work shall also be devoted to the application of Active Flux to other nonlinear balance laws with more sophisticated stationary states.

\section*{Acknowledgement}

WB was supported by the Deutsche Forschungsgemeinschaft (DFG) through project 429491391 (BA 6878/1-1). JB was supported by the Klaus Tschira Foundation. 

\bibliographystyle{alpha}

\begin{thebibliography}{CMDG{\etalchar{+}}12}

\bibitem[ABB{\etalchar{+}}04]{Audusse2004}
Emmanuel Audusse, Fran{\c{c}}ois Bouchut, Marie-Odile Bristeau, Rupert Klein,
  and Benoit Perthame.
\newblock A fast and stable well-balanced scheme with hydrostatic
  reconstruction for shallow water flows.
\newblock {\em SIAM Journal on Scientific Computing}, 25(6):2050--2065, 2004.

\bibitem[Bar21]{barsukow19activeflux}
Wasilij Barsukow.
\newblock The active flux scheme for nonlinear problems.
\newblock {\em Journal of Scientific Computing}, 86(1):1--34, 2021.

\bibitem[BBK21]{barsukow19activefluxsource}
Wasilij Barsukow, Jonas~P Berberich, and Christian Klingenberg.
\newblock On the {A}ctive {F}lux scheme for hyperbolic {PDE}s with source
  terms.
\newblock {\em SIAM Journal on Scientific Computing}, 43(6):A4015--A4042, 2021.

\bibitem[BCK19]{Berberich2019b}
Jonas~P Berberich, Praveen Chandrashekar, and Christian Klingenberg.
\newblock High order well-balanced finite volume methods for multi-dimensional
  systems of hyperbolic balance laws.
\newblock {\em arXiv preprint arXiv:1903.05154}, 2019.

\bibitem[BCKN13]{bollermann13}
Andreas Bollermann, Guoxian Chen, Alexander Kurganov, and Sebastian Noelle.
\newblock A well-balanced reconstruction of wet/dry fronts for the shallow
  water equations.
\newblock {\em Journal of Scientific Computing}, 56(2):267--290, 2013.

\bibitem[BdSV71]{saintvenant71}
Adh{\'e}mar Jean-Claude Barré~de Saint-Venant.
\newblock Th{\'e}orie du mouvement non-permanent des eaux, avec application aux
  crues des rivières et á l’introduction des mar{\'e}es dans leur lit.
\newblock {\em Comptes Rendus Acad. Sci. Paris}, 73(147-154):8, 1871.

\bibitem[BHKR19]{barsukow18activeflux}
Wasilij Barsukow, Jonathan Hohm, Christian Klingenberg, and Philip~L Roe.
\newblock The active flux scheme on {C}artesian grids and its low {M}ach number
  limit.
\newblock {\em Journal of Scientific Computing}, 81(1):594--622, 2019.

\bibitem[BK22]{barsukow17}
Wasilij Barsukow and Christian Klingenberg.
\newblock Exact solution and a truly multidimensional {G}odunov scheme for the
  acoustic equations.
\newblock {\em ESAIM: M2AN}, 56(1), 2022.

\bibitem[BNLM11]{bollermann11}
Andreas Bollermann, Sebastian Noelle, and Maria
  Luk{\'a}{\v{c}}ov{\'a}-Medvid’ov{\'a}.
\newblock Finite volume evolution {G}alerkin methods for the shallow water
  equations with dry beds.
\newblock {\em Communications in Computational Physics}, 10(2):371--404, 2011.

\bibitem[Bou04]{bouchut04}
Fran{\c{c}}ois Bouchut.
\newblock {\em Nonlinear stability of finite Volume Methods for hyperbolic
  conservation laws and Well-Balanced schemes for sources}.
\newblock Springer Science \& Business Media, 2004.

\bibitem[BV94]{Bermudez1994}
Alfredo Bermudez and Elena V{\'a}zquez.
\newblock Upwind methods for hyperbolic conservation laws with source terms.
\newblock {\em Computers \& Fluids}, 23(8):1049--1071, 1994.

\bibitem[CCH{\etalchar{+}}19]{cheng19}
Yuanzhen Cheng, Alina Chertock, Michael Herty, Alexander Kurganov, and Tong Wu.
\newblock A new approach for designing moving-water equilibria preserving
  schemes for the shallow water equations.
\newblock {\em Journal of Scientific Computing}, 80(1):538--554, 2019.

\bibitem[CGVP06]{castro06}
Manuel~J Castro, Jos{\'e}~M Gonz{\'a}lez-Vida, and Carlos Par{\'e}s.
\newblock Numerical treatment of wet/dry fronts in shallow flows with a
  modified {R}oe scheme.
\newblock {\em Mathematical Models and Methods in Applied Sciences},
  16(06):897--931, 2006.

\bibitem[CLS04]{chinnayya04}
Ashwin Chinnayya, Alain-Yves LeRoux, and Nicolas Seguin.
\newblock A well-balanced numerical scheme for the approximation of the
  shallow-water equations with topography: the resonance phenomenon.
\newblock {\em Int. J. Finite}, 1(1):33, 2004.

\bibitem[CMDG{\etalchar{+}}12]{cozzolino12}
Luca Cozzolino, Renata~Della Morte, Giuseppe Del~Giudice, Anna Palumbo, and
  Domenico Pianese.
\newblock A well-balanced spectral volume scheme with the wetting--drying
  property for the shallow-water equations.
\newblock {\em Journal of Hydroinformatics}, 14(3):745--760, 2012.

\bibitem[CS19]{castro19}
Manuel~J Castro and Matteo Semplice.
\newblock Third-and fourth-order well-balanced schemes for the shallow water
  equations based on the {CWENO} reconstruction.
\newblock {\em International Journal for Numerical Methods in Fluids},
  89(8):304--325, 2019.

\bibitem[CWY03]{carrier03}
George~F Carrier, Tai~Tei Wu, and Harry Yeh.
\newblock Tsunami run-up and draw-down on a plane beach.
\newblock {\em Journal of Fluid Mechanics}, 475:79, 2003.

\bibitem[ER11]{eymann11}
Timothy~A Eymann and Philip~L Roe.
\newblock Active flux schemes for systems.
\newblock In {\em 20th AIAA computational fluid dynamics conference}, 2011.

\bibitem[ER13]{eymann13}
Timothy~A Eymann and Philip~L Roe.
\newblock Multidimensional active flux schemes.
\newblock In {\em 21st AIAA computational fluid dynamics conference}, 2013.

\bibitem[Eym13]{eymann13a}
Timothy~Andrew Eymann.
\newblock Active flux schemes.
\newblock {\em PhD thesis, University of Michigan}, 2013.

\bibitem[Fan17]{fan17}
Duoming Fan.
\newblock {\em On the acoustic component of active flux schemes for nonlinear
  hyperbolic conservation laws}.
\newblock PhD thesis, University of Michigan, Dissertation, 2017.

\bibitem[GHS03]{gallouet03}
Thierry Gallou{\"e}t, Jean-Marc H{\'e}rard, and Nicolas Seguin.
\newblock Some approximate {G}odunov schemes to compute shallow-water equations
  with topography.
\newblock {\em Computers \& Fluids}, 32(4):479--513, 2003.

\bibitem[Gos01]{gosse01}
Laurent Gosse.
\newblock A well-balanced scheme using non-conservative products designed for
  hyperbolic systems of conservation laws with source terms.
\newblock {\em Mathematical Models and Methods in Applied Sciences},
  11(02):339--365, 2001.

\bibitem[GPC07]{gallardo07}
Jos{\'e}~M Gallardo, Carlos Par{\'e}s, and Manuel Castro.
\newblock On a well-balanced high-order finite volume scheme for shallow water
  equations with topography and dry areas.
\newblock {\em Journal of Computational Physics}, 227(1):574--601, 2007.

\bibitem[HKS19]{kerkmann18}
Christiane Helzel, David Kerkmann, and Leonardo Scandurra.
\newblock A new {ADER} method inspired by the active flux method.
\newblock {\em Journal of Scientific Computing}, 80(3):1463--1497, 2019.

\bibitem[KL02]{kurganov02}
Alexander Kurganov and Doron Levy.
\newblock Central-upwind schemes for the {S}aint-{V}enant system.
\newblock {\em ESAIM: Mathematical Modelling and Numerical Analysis},
  36(3):397--425, 2002.

\bibitem[Kur18]{kurganov18}
Alexander Kurganov.
\newblock Finite-volume schemes for shallow-water equations.
\newblock {\em Acta Numerica}, 27:289--351, 2018.

\bibitem[LeV98]{leveque98}
Randall~J LeVeque.
\newblock Balancing source terms and flux gradients in high-resolution
  {G}odunov methods: the quasi-steady wave-propagation algorithm.
\newblock {\em Journal of computational physics}, 146(1):346--365, 1998.

\bibitem[PC04]{pares04}
Carlos Par{\'e}s and Manuel Castro.
\newblock On the well-balance property of {R}oe's method for nonconservative
  hyperbolic systems. {A}pplications to shallow-water systems.
\newblock {\em ESAIM: mathematical modelling and numerical analysis},
  38(5):821--852, 2004.

\bibitem[PR17]{pareschi2017}
Lorenzo Pareschi and Thomas Rey.
\newblock Residual equilibrium schemes for time dependent partial differential
  equations.
\newblock {\em Computers \& Fluids}, 156:329--342, 2017.

\bibitem[RB09]{ricchiuto09}
Mario Ricchiuto and Andreas Bollermann.
\newblock Stabilized residual distribution for shallow water simulations.
\newblock {\em Journal of Computational Physics}, 228(4):1071--1115, 2009.

\bibitem[RBT03]{rogers03}
Benedict~D Rogers, Alistair~GL Borthwick, and Paul~H Taylor.
\newblock Mathematical balancing of flux gradient and source terms prior to
  using {R}oe’s approximate {R}iemann solver.
\newblock {\em Journal of Computational Physics}, 192(2):422--451, 2003.

\bibitem[RLM15]{roe15}
Philip~L Roe, Tyler Lung, and Jungyeoul Maeng.
\newblock New approaches to limiting.
\newblock In {\em 22nd AIAA Computational Fluid Dynamics Conference}, page
  2913, 2015.

\bibitem[Roe87]{roe87}
PL~Roe.
\newblock Upwind differencing schemes for hyperbolic conservation laws with
  source terms.
\newblock In {\em Nonlinear hyperbolic problems}, pages 41--51. Springer, 1987.

\bibitem[Roe21]{roe21}
Philip Roe.
\newblock Designing cfd methods for bandwidth—a physical approach.
\newblock {\em Computers \& Fluids}, 214:104774, 2021.

\bibitem[Tha81]{thacker81}
William~Carlisle Thacker.
\newblock Some exact solutions to the nonlinear shallow-water wave equations.
\newblock {\em Journal of Fluid Mechanics}, 107:499--508, 1981.

\bibitem[VB14]{vater14}
Stefan Vater and J{\"o}rn Behrens.
\newblock Well-balanced inundation modeling for shallow-water flows with
  discontinuous {G}alerkin schemes.
\newblock In {\em Finite volumes for complex applications VII-elliptic,
  parabolic and hyperbolic problems}, pages 965--973. Springer, 2014.

\bibitem[vL77]{vanleer77}
Bram van Leer.
\newblock Towards the ultimate conservative difference scheme. {IV}. {A} new
  approach to numerical convection.
\newblock {\em Journal of computational physics}, 23(3):276--299, 1977.

\bibitem[XS05]{xing05}
Yulong Xing and Chi-Wang Shu.
\newblock High order finite difference {WENO} schemes with the exact
  conservation property for the shallow water equations.
\newblock {\em Journal of Computational Physics}, 208(1):206--227, 2005.

\bibitem[XS06]{xing06}
Yulong Xing and Chi-Wang Shu.
\newblock High order well-balanced finite volume {WENO} schemes and
  discontinuous {G}alerkin methods for a class of hyperbolic systems with
  source terms.
\newblock {\em Journal of Computational Physics}, 214(2):567--598, 2006.

\bibitem[XS11]{xing11}
Yulong Xing and Chi-Wang Shu.
\newblock High-order finite volume {WENO} schemes for the shallow water
  equations with dry states.
\newblock {\em Advances in Water Resources}, 34(8):1026--1038, 2011.

\bibitem[XS14]{xing14}
Yulong Xing and Chi-Wang Shu.
\newblock A survey of high order schemes for the shallow water equations.
\newblock {\em J. Math. Study}, 47(3):221--249, 2014.

\bibitem[Xu02]{xu02}
Kun Xu.
\newblock A well-balanced gas-kinetic scheme for the shallow-water equations
  with source terms.
\newblock {\em Journal of Computational Physics}, 178(2):533--562, 2002.

\bibitem[XZS10]{xing10}
Yulong Xing, Xiangxiong Zhang, and Chi-Wang Shu.
\newblock Positivity-preserving high order well-balanced discontinuous
  {G}alerkin methods for the shallow water equations.
\newblock {\em Advances in Water Resources}, 33(12):1476--1493, 2010.

\bibitem[ZCMI01]{zhou01}
Jian~G Zhou, Derek~M Causon, Clive~G Mingham, and David~M Ingram.
\newblock The surface gradient method for the treatment of source terms in the
  shallow-water equations.
\newblock {\em Journal of Computational physics}, 168(1):1--25, 2001.

\bibitem[ZS10]{zhang10}
Xiangxiong Zhang and Chi-Wang Shu.
\newblock On positivity-preserving high order discontinuous {G}alerkin schemes
  for compressible {E}uler equations on rectangular meshes.
\newblock {\em Journal of Computational Physics}, 229(23):8918--8934, 2010.

\end{thebibliography}
\newcommand{\etalchar}[1]{$^{#1}$}

\appendix

\section{Details on the reconstruction procedure} \label{app:reconstruction}

\newcommand{\eqn}{\mathscr{G}}

In the following, define a relative coordinate in the cell
 \begin{align}
  \xi := \frac12 - \frac{x^*}{\Delta x} \in (0,1)
 \end{align}
 Denote by $b_i$ the coefficients in $b(x) = b_0 + b_1 x + b_2 x^2$ and define the two parameters
 \begin{align}
  \bar h_\text{crit} &:= h_\text L + \sqrt{\frac{(h_\text R - h_\text L)^3}{9 \lvert b_2\rvert \Delta x^2 }} &
  \xi_0 &:= \sqrt{\frac{h_\text R - h_\text L}{\|b_2\| \Delta x^2} }
 \end{align}

\begin{figure}
 \centering
 \includegraphics[width=\textwidth]{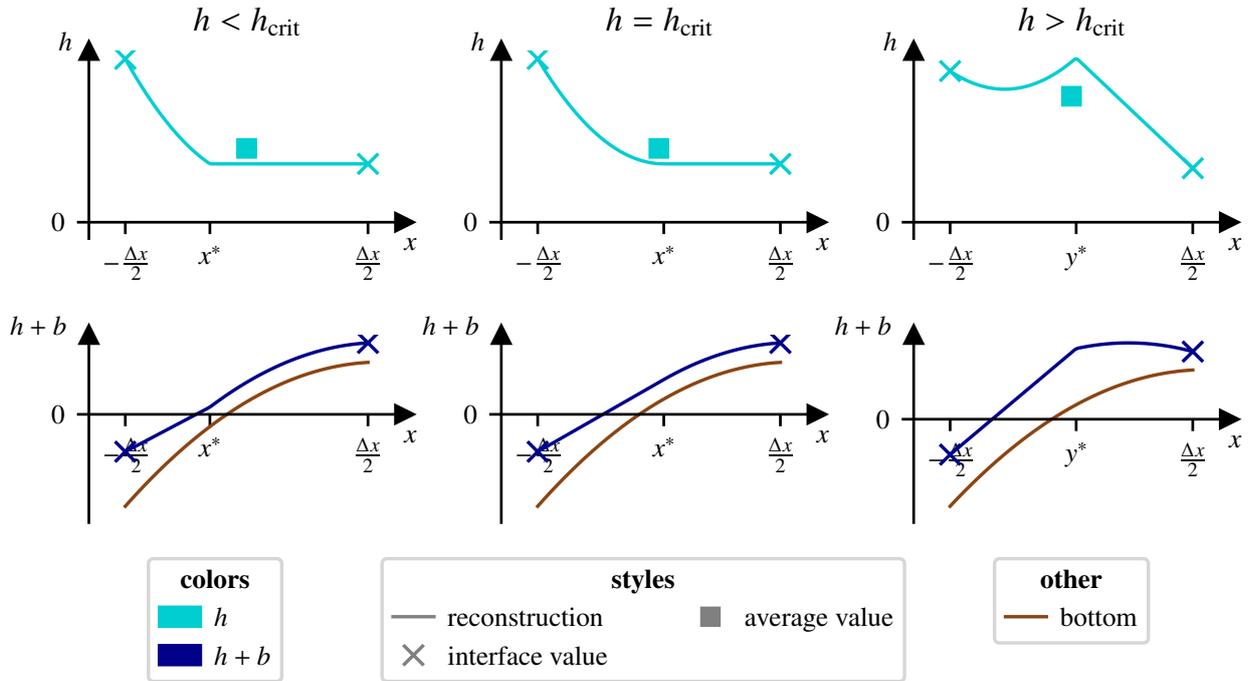}
 \caption{Reconstruction according Lemma \ref{lem:recondrywet}. \emph{Top}: Water surface ($h+b$) is shown. \emph{Bottom}: Water height $h$ is shown. $\bar h$ varies from $\bar h =0.19$ (\emph{left}, reconstruction \eqref{eq:reconwetdrygeneric}), over $\bar h = \bar h_\text{crit} \simeq 0.196765$ (\emph{center}), to $\bar h = 0.35$ (\emph{right}, reconstruction \eqref{eq:reconwetdryexceptional} because no $x^*$ for a reconstruction according to \eqref{eq:reconwetdrygeneric} exists). Clearly, a parabolic reconstruction is given preference whenever it is positive, while reconstructions \eqref{eq:reconwetdrygeneric}--\eqref{eq:reconwetdryexceptional} cover the remaining cases. In this figure, however, to improve visibility, the reconstruction \eqref{eq:reconwetdryexceptional} (bottom row) is shown for a choice of the average and of point values which would actually allow for a well-behaved parabolic reconstruction. The actual regime of application of \eqref{eq:reconwetdryexceptional} occurs in a small parameter range that makes visualization difficult. The interested reader finds it depicted accurately, but much less distinctly in Figure \ref{fig:casebrightsmallconcave}.}
 \label{fig:recondrywetLeft}
\end{figure}

\begin{lemma} \label{lem:recondrywet}

 Assume $h_\text L < \bar h < h_\text R$ and assume also that the parabolic reconstruction is violating non-negativity. If $b_2 < 0$, $\xi_0 < 1$ and $\bar h > \bar h_\text{crit}$, then the reconstruction
 \begin{align}
  h_\text{recon}(x) = \begin{cases} h_\text L + (x + \frac{\Delta x}{2}) \frac{h^* - h_\text L}{y^* + \frac{\Delta x}{2}} & -\frac{\Delta x}{2} \leq x < y^*\\
                       h^* + b(y^*) + (x - y^*) \frac{h_\text R + b_\text R - h^* - b(y^*)}{\frac{\Delta x}{2} - y^*} - b(x) & y^* \leq x < \frac{\Delta x}{2}
                      \end{cases} \label{eq:reconwetdryexceptional}
 \end{align}
 with $y^* = \frac{\Delta x}2 - \Delta x \xi_0$ and $h^* = 2 (\bar h - \bar h_\text{crit} )+ h_\text L$ is conservative and positive, and otherwise there exists $x^* \in \left(-\frac{\Delta x}{2},\frac{\Delta x}{2}\right)$ such that the reconstruction
 \begin{align}
  h_\text{recon}(x) = \begin{cases} h_\text L & -\frac{\Delta x}{2} \leq x < x^*\\
                       h_\text L + b(x^*) + (x - x^*) \frac{h_\text R + b_\text R - h_\text L - b(x^*)}{\frac{\Delta x}{2} - x^*} - b(x) & x^* \leq x < \frac{\Delta x}{2}
                      \end{cases} \label{eq:reconwetdrygeneric}
 \end{align}
 is conservative, positive and well-balanced at the dry/wet interface (see Figure \ref{fig:recondrywetLeft}).
\end{lemma} 

\emph{Note}: In the case of a lake at rest, reconstruction \eqref{eq:reconwetdrygeneric} is exact, and thus exists. 

\begin{proof} 
 The proof starts with reconstruction \eqref{eq:reconwetdrygeneric}, and shows how the exceptional case \eqref{eq:reconwetdryexceptional} arises. Both reconstructions obviously fulfill $h_\text{recon}\left( \pm \frac{\Delta x}{2} \right) = h_\text{L/R}$. The reconstruction \eqref{eq:reconwetdrygeneric} is conservative if $\xi$ fulfills
 \begin{align}
  0 &=  -(\bar h - h_\text L) + \frac{h_\text R - h_\text L}{2}\xi + \frac16 b_2 \Delta x^2 \xi^3 =: \eqn(\xi) \label{eq:xilocation}
 \end{align}
 which follows from
 \begin{align}
  \Delta x \bar h &=  \int^{\frac{\Delta x}{2}}_{-\frac{\Delta x}{2}} h_\text{recon}(x) \, \dd x = 
  \Delta x h_\text L + \frac{h_\text R - h_\text L}{2} \left(\frac{\Delta x}{2} - x^*\right) + \frac16 b_2 \left(\frac{\Delta x}{2} - x^*\right)^3
 \end{align}
 This is a cubic equation in $\xi$ (a so-called \emph{depressed} cubic, even). The existence of a suitable real solution is established next.

 It is clear, both geometrically and analytically, that as $\bar h \to h_\text L$, $x^* \to \frac{\Delta x}{2}$, i.e.\ $\xi \to 0$. Computing the discriminant, equation \ref{eq:xilocation} is found to have one real solution if $-  \left( \frac{h_\text R - h_\text L}{b_2 \Delta x^2 } \right)^3 - 9 \left(  \frac{h_\text L - \bar h}{b_2 \Delta x^2 } \right )^2 < 0$, two in case of equality and three otherwise. The equation thus has a unique real solution if
 \begin{align}
  -   \frac{(h_\text R - h_\text L)^3}{9 b_2 \Delta x^2 } <  (\bar h - h_\text L)^2 \label{eq:discriminantsimple}
 \end{align}
 Next, two cases need to be considered, according to the sign of $b_2$.
 \begin{enumerate}[I.]
 \item If $b_2$ is positive, then \eqref{eq:discriminantsimple} is always true. Note that $\eqn(0) < 0$, and for $b_2 > 0$, $\eqn(\xi) \to +\infty$ as $\xi \to +\infty$. Thus, the unique real solution of \eqref{eq:xilocation} is located in $(0, \infty)$ (see Figure \ref{fig:Greconallcases}, \emph{left}). Also, the solution of $\eqn(\xi)=0$ is a monotonically increasing function of $\bar h$, which can be seen by differentiating \eqref{eq:xilocation} with respect to $\bar h$:
 \begin{align}
  1 &= \frac12 \left( h_\text R - h_\text L  + b_2 \Delta x^2 \xi^2 \right ) \frac{\dd \xi}{\dd \bar h} \label{eq:Gmonotone}
 \end{align}
 Starting from $h_\text L$, as $\bar h$ increases, $\xi$ is monotonically increasing until $x^* = - \frac{\Delta x}{2}$ ($\xi = 1$) which occurs when
 \begin{align}
  \bar h = \frac{h_\text R + h_\text L}{2} + \frac16 b_2 \Delta x^2 > \frac{h_\text R + h_\text L}{2}
 \end{align}
 This case will never actually be reached, because the parabolic reconstruction is monotone and does not violate non-negativity for even smaller values of $\bar h$. It will thus be given preference.
 
 Finally, the parabolic reconstruction in $(x^*, \frac{\Delta x}{2})$ is positive at the endpoints of the interval, and concave. It thus is positive in the entire interval. This proves the statement of the Lemma for $b_2 > 0$.

 \item If $b_2$ is negative, then there is a unique solution only if
 \begin{align}
   \bar h > h_\text L + \sqrt{\frac{(h_\text R - h_\text L)^3}{9 \lvert b_2\rvert \Delta x^2 }} =: \bar h_\text{crit}
 \end{align}
 Still $\eqn(0) < 0$, but now, for $b_2 < 0$, $\eqn(\xi) \to +\infty$ as $\xi \to -\infty$. Thus, there is always one real solution of \eqref{eq:xilocation} located in $(-\infty, 0)$. This solution, however, cannot be used because $\xi$ is required to be in $(0,1)$. This, again, requires distinguishing two cases.
 \begin{enumerate}[a.]
 \item If $\bar h < \bar h_\text{crit}$, then there are two real solutions of \eqref{eq:xilocation} located in $\mathbb R^+$ (see Figure \ref{fig:Greconallcases}, center), and the smallest can be used for the reconstruction. 
 By differentiating $\eqn(\xi)$ with respect to $\xi$ it can easily be seen that this solution is situated in $(0, \xi_0)$. Thus, when it is chosen,
 \begin{align}
  \xi^2 < \frac{h_\text R - h_\text L}{\lvert b_2\rvert \Delta x^2}
 \end{align} 
 and from \eqref{eq:Gmonotone} one finds that again, this solution of $\eqn(\xi)=0$ is a monotonically increasing function of $\bar h$. Starting from $h_\text L$, as $\bar h$ increases, $\xi$ increases until either $\xi = 1$, i.e.
 \begin{align}
 \bar h = \frac{h_\text R + h_\text L}{2} - \frac16 \lvert b_2\rvert \Delta x^2 \label{eq:GxiOneb2neg}
 \end{align}
or $\bar h = \bar h_\text{crit}$, whatever happens first.
 

 \begin{figure}
  \centering
  \includegraphics[width=\textwidth]{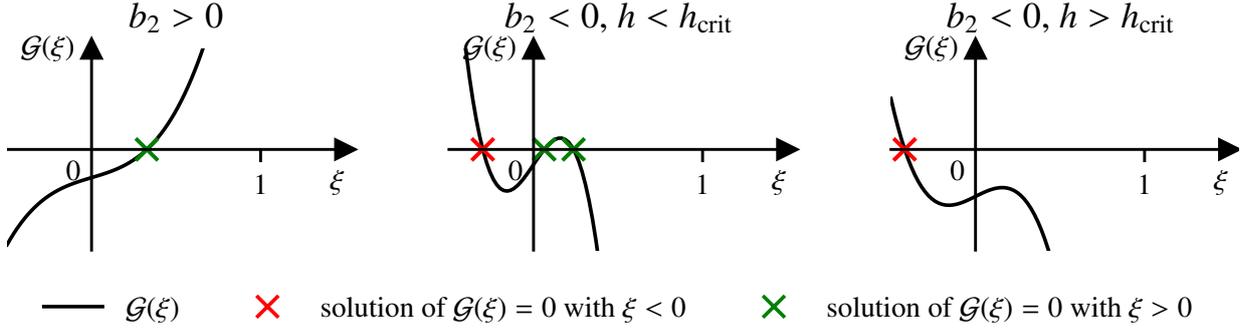}
  \caption{$\eqn(\xi)$ shown for representative combinations of parameters. \emph{Left}: $b_2 > 0$, $\eqn(\xi)=0$ has a unique solution, located in $(0, \infty)$. \emph{Center}: $b_2 < 0$, $\bar h < \bar h_\text{crit}$, $\eqn(\xi)=0$ has two positive solutions. \emph{Right}: $b_2 < 0$, $\bar h > \bar h_\text{crit}$, and $\eqn(\xi) =0$ has no positive solutions.}
  \label{fig:Greconallcases}
 \end{figure}

 If $\xi = 1$, then the reconstruction inside the cell is linear and necessarily equal to the parabolic reconstruction. It thus does not violate non-negativity. This can also be shown explicitly, by noticing that \eqref{eq:GxiOneb2neg} and $\xi < \xi_0$ together violate the second condition in \eqref{eq:reconnonneg}.
 
 The parabolic part of the reconstruction is convex, and the value at the (unique) minimum of the parabola is located at
 \begin{align}
  -\frac{(h_\text R - h_\text L) - \lvert b_2\rvert  \Delta x^2 \xi (1 - \xi)}{2 \lvert b_2\rvert \Delta x \xi} < - \frac{\xi \Delta x}{2} < (1-2\xi) \frac{\Delta x}{2}
 \end{align}
 where $\xi < \xi_0$ was used in the first inequality and $- \xi < 1 - 2 \xi$ for $\xi \in (0,1)$ in the second. Thus, the parabola is monotone inside $(x^*, \frac{\Delta x}{2})$, and positive.
 
 \item If $\bar h >  \bar h_\text{crit}$, the real solution to \eqref{eq:xilocation} is unique (see Figure \ref{fig:Greconallcases}, right), but located outside $(0,1)$ and the reconstruction procedure is not possible. If $\bar h = \bar h_\text{crit}$, then the relevant solution $\xi$ of $\eqn(\xi) = 0$ is also a solution of $\eqn'(\xi) = 0$ which is $\xi = \xi_0$. Then, the reconstruction is $C^1$ at $x^*$, because
 \begin{align}
  h_\text{recon}'(x^*) = \frac{h_\text R + b_\text R - h_\text L - b(x^*)}{\frac{\Delta x}{2} - x^*} - b'(x^*) = 0
 \end{align}

 If $\bar h >  \bar h_\text{crit}$, the reconstruction \eqref{eq:reconwetdrygeneric} cannot always be used. It requires a modification if $\xi_0 < 1$ and $\bar h_\text{crit}$ fulfills the conditions of Lemma \ref{lem:nonnegativityparabolic}. From $\frac{h_\text R - h_\text L}{\lvert b_2 \rvert \Delta x^2} < 1$ follows
 \begin{align}
  (\bar h_\text{crit} - h_\text L)^2 &= \frac{(h_\text R - h_\text L)^3}{9\lvert b_2 \rvert \Delta x^2} < \left(\frac{h_\text R - h_\text L}{3} \right )^2\\
  \bar h_\text{crit}  <  \frac{h_\text L + h_\text R }{2} - \frac{h_\text R - h_\text L}{6}
 \end{align}
 which is the second condition of Lemma \ref{lem:nonnegativityparabolic} (equation \eqref{eq:reconnonneg}). It thus is always fulfilled in this case. Instead of evaluating exactly the conditions under which the first condition of Lemma \ref{lem:nonnegativityparabolic} is fulfilled as well, we suggest in this case to generally replace the reconstruction \eqref{eq:reconwetdrygeneric} by the reconstruction \eqref{eq:reconwetdryexceptional}:

 \begin{align*}
  h_\text{recon}(x) = \begin{cases} h_\text L + (x + \frac{\Delta x}{2}) \frac{h^* - h_\text L}{y^* + \frac{\Delta x}{2}} & -\frac{\Delta x}{2} \leq x < y^*\\
                       h^* + b(y^*) + (x - y^*) \frac{h_\text R + b_\text R - h^* - b(y^*)}{\frac{\Delta x}{2} - y^*} - b(x) & y^* \leq x < \frac{\Delta x}{2}
                      \end{cases}
 \end{align*}
 where now $y^* = \frac{\Delta x}2 - \Delta x \xi_0$ and 
 \begin{align}
  h^* = 2 \bar h - h_\text L - 2\sqrt{ \frac{(h_\text R - h_\text L)^3}{9 \lvert b_2 \rvert \Delta x^2}} = 2 (\bar h - \bar h_\text{crit} )+ h_\text L  > h_\text L
 \end{align}
 which ensures that this reconstruction is conservative. As $h^* > h_\text L$, the new reconstruction is positive as well.

 \end{enumerate}
 \end{enumerate}
 Note that, given the above analysis, actually solving \eqref{eq:xilocation} in the real numbers is trivial using Cardano's and Vieta's formulae.
\end{proof}

The mirror case, i.e.\ when $h_\text L < \bar h < h_\text R$, is obtained analogously. This time, it is convenient to define $\xi := \frac12 + \frac{x^*}{\Delta x}$ and
\begin{align}
 \bar h_\text{crit} &:= h_\text R + \sqrt{\frac{(h_\text L - h_\text R)^3}{9 \lvert b_2 \rvert \Delta x^2}} &
 \xi_0 &:= \sqrt{\frac{h_\text L - h_\text R}{\lvert b_2 \rvert \Delta x^2}}
\end{align}

\begin{lemma} \label{lem:recondrywetmirror}
If $b_2 < 0$, $\bar h > h_\text{crit}$ and $\xi_0 < 1$, then the reconstruction is
\begin{align}
 h_\text{recon}(x) &= \begin{cases}
                      h_\text L + b_\text L + (x + \frac{\Delta x}{2}) \frac{h^* + b(y^*) - h_\text L - b_\text L}{y^* + \frac{\Delta x}{2}} - b(x) & -\frac{\Delta x}{2} \leq x < y^*\\
                      h^* + (x - y^*) \frac{h_\text R - h^*}{\frac{\Delta x}{2} - y^*} & y^* \leq x < \frac{\Delta x}{2}
                     \end{cases} \label{eq:reconwetdryexceptionalmirror}
\end{align}
with $y^* = \Delta x (\xi_0 - \frac12)$ and $h^* = 2 (\bar h - \bar h_\text{crit}) + h_\text R$, while otherwise there exists $x^* \in (-\frac{\Delta x}{2}, \frac{\Delta x}{2})$
\begin{align}
 h_\text{recon}(x) &= \begin{cases}
                      h_\text L + b_\text L + (x + \frac{\Delta x}{2}) \frac{h_\text R + b(x^*) - h_\text L - b_\text L}{x^* + \frac{\Delta x}{2}} - b(x) & -\frac{\Delta x}{2} \leq x < x^*\\
                      h_\text R & x^* \leq x < \frac{\Delta x}{2}
                     \end{cases} \label{eq:reconwetdrygenericmirror}
\end{align}
\end{lemma}
\begin{proof}
 The proof is entirely analogous to that of Lemma \ref{lem:recondrywet}. Again, if $b_2 > 0$, there is a unique real solution that can be used for the reconstruction. If $b_2 < 0$, as long as there are three real solutions $\xi$, one them is located in $(0, \xi_0)$ and can be used. Only one real root in this case exists if $\bar h > h_\text{crit}$. This, however only requires action, if the location $\xi_0 < 1$. This justifies \eqref{eq:reconwetdryexceptionalmirror}.
\end{proof}

 \begin{lemma} \label{lem:positivewaterheight}
 The non-negative reconstruction has the following properties:
  \begin{enumerate}[i.]
   \item \label{it:lemx1} For case \ref{it:lowerthanboth}: $\displaystyle - \frac{\Delta x}{2} <  x^*_1 < 0$.
   \item \label{it:lemmosesmom} The reconstruction \eqref{eq:mosesmomentum} of the momentum is conservative.
   \item \label{it:lemhalf1} The reconstruction \eqref{eq:recononeflatx2} and \eqref{eq:recononeflatx3} of the momentum is conservative.
  \end{enumerate}
 \end{lemma}

\textit{Proof}.
  \begin{enumerate}[i.]
   \item Using $\bar h > 0$ and $2 \bar h < h_\text L + h_\text R$:
  \begin{align}
   x^*_1  > - \frac{\Delta x}{2}  \frac{ 1}{1 - \frac{2 \bar h}{h_\text L + h_\text R} f  }  > - \frac{\Delta x}{2} 
  \end{align}
  Observe also that because of $0 < f < 1$ one has $\bar h f < \bar h < \frac{h_\text L + h_\text R}{2}$ and
  \begin{align}
     x^*_1 = -\Delta x  \frac{ 2 \bar h  \left( f - 2 + \frac{h_\text L + h_\text R}{2 \bar h} \right )}{4 \left( \frac{h_\text L + h_\text R}{2} -\bar h f \right )} \leq 0
  \end{align}
  if $f - 2 + \frac{h_\text L + h_\text R}{2 \bar h} \geq 0$.

  \item Integrate \eqref{eq:mosesmomentum} over the cell $[-\frac{\Delta x}{2}, \frac{\Delta x}{2}]$ and use the property \eqref{eq:reconconservation} of $\mathscr R$.
  
   \begin{align}
    \int_{-\frac{\Delta x}{2}}^{\frac{\Delta x}{2}} \dd x m_\text{recon}(x) &= 
    \int_{-\frac{\Delta x}{2}}^{x^*_\text L} \mathscr R\left( \bar m, m_\text L, \bar m; -\frac{\Delta x}{2}, x^*_\text L\right)  \dd x 
    \\&+ \bar m (x^*_\text R - x^*_\text L)
    + \int_{x^*_\text R}^{\frac{\Delta x}{2}}\mathscr R\left( \bar m, \bar m, m_\text R;  x^*_\text R, \frac{\Delta x}{2}; x\right) \dd x\\
    &=  \bar m \left( x^*_\text L + \frac{\Delta x}{2}\right ) + \bar m (x^*_\text R - x^*_\text L) +  \bar m \left( \frac{\Delta x}{2} - x^*_\text R \right )\\
    &=  \bar m \left( x^*_\text L  - x^*_\text R + \Delta x \right )  + \bar m (x^*_\text R - x^*_\text L) = \bar m \Delta x
   \end{align}
   
  \item Integrate the reconstruction \eqref{eq:recononeflatx2} over $[-\frac{\Delta x}{2}, \frac{\Delta x}{2}]$:
  
   \begin{align} 
    \Delta x \bar m &= \int_{-\frac{\Delta x}{2}}^{\frac{\Delta x}{2}} \dd x \, m_\text{recon}(x) \\
   &= \int_{-\frac{\Delta x}{2}}^{x^*_\text L} \dd x \, 
	\mathscr R\left(\frac{ \bar m  - (\frac{1}{2} - \frac{x^*_\text L}{\Delta x}) m_\text R}{\frac{1}{2} + \frac{x^*_\text L}{\Delta x}},  m_\text L,m_\text R; -\frac{\Delta x}{2}, x^*_\text L; x\right) + m_\text R (\frac{\Delta x}{2} -x^*_\text L)\\
 	&= \frac{ \bar m \Delta x - (\frac{\Delta x}{2} - x^*_\text L) m_\text R}{\frac{\Delta x}{2} + x^*_\text L} \left( \frac{\Delta x}{2} + x^*_\text L \right ) + m_\text R (\frac{\Delta x}{2} -x^*_\text L)
    \end{align}
    
    Analogously, upon integrating \eqref{eq:recononeflatx3}

   \begin{align}
    \Delta x \bar m &= \int_{-\frac{\Delta x}{2}}^{\frac{\Delta x}{2}} \dd x \, m_\text{recon}(x) \\
    &= m_\text L (x^*_\text R + \frac{\Delta x}{2}) + \frac{\bar m \Delta x -  m_\text L (x^*_\text R + \frac{\Delta x}{2})}{\frac{\Delta x}{2} - x^*_\text R}  (\frac{\Delta x}{2} - x^*_\text R)
   \end{align}

  \end{enumerate}

\vfill

\phantom{m}

\includegraphics[width=\textwidth]{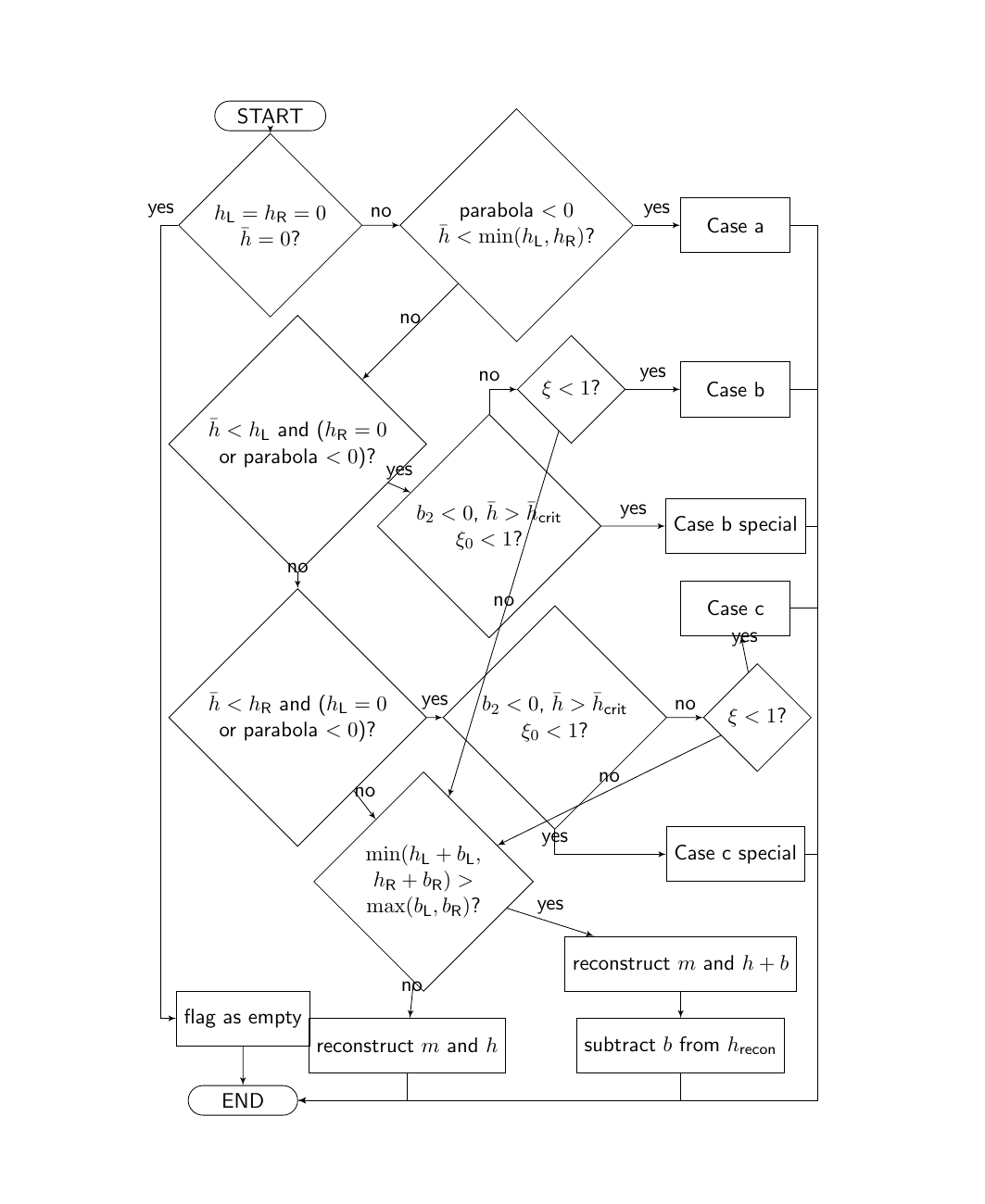}

\newpage

\section{Details on the point value update}

The following flow chart summarizes all the steps of the point value update.

\vspace{0.5cm}

\includegraphics[width=\textwidth]{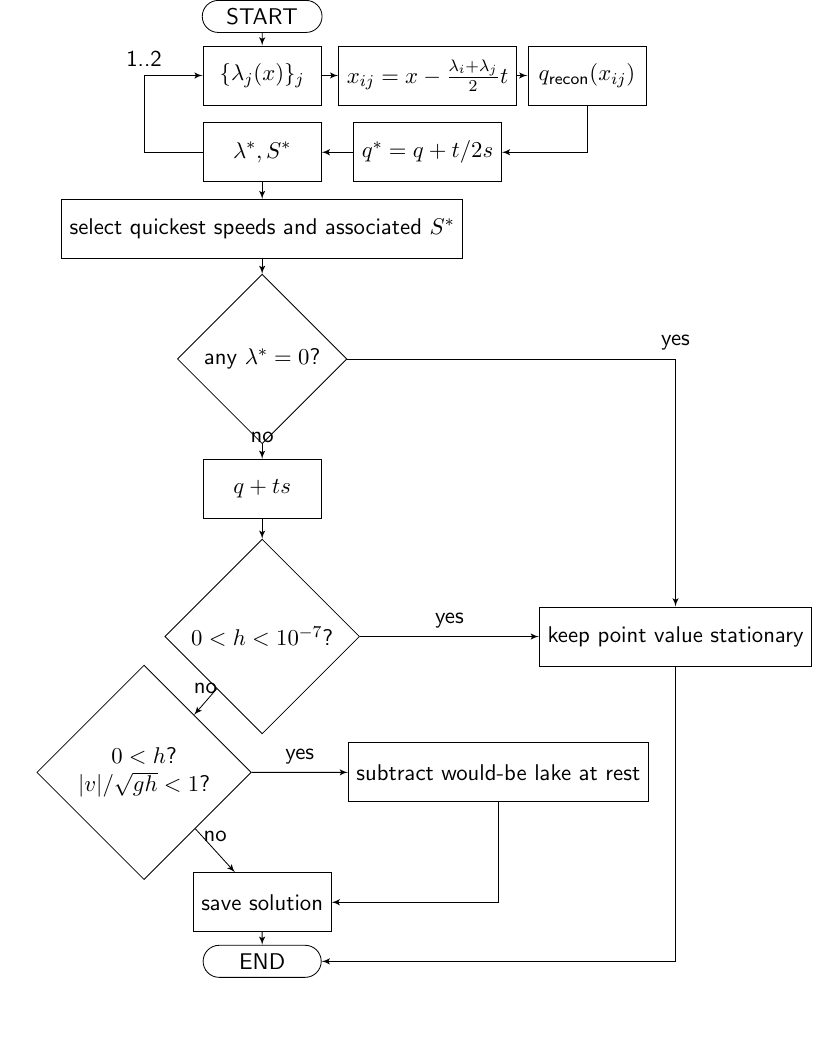}

\newpage

\section{Regularization of $h=0$} \label{ssec:drycellsregularizationpointvalues}

A practical aspect of dealing with dry states is the necessity to replace all expressions which are singular at $h=0$ by some regularized expressions. A numerical code might otherwise crash. These corrections do not need to be strictly physical. Therefore, if $h\leq 0$, then
\begin{itemize}
 \item The characteristic variables are set to 0.
 \item The characteristic speeds are set to 0 (also if\footnote{The rationale is that otherwise we will be dividing by a small number (which probably has arisen from random errors) to obtain a speed, and we have observed these speeds to be completely unphysical.} $h < 10^{-14}$).
 \item The sound speed and pressure is set to 0.
 \item The flux $f(q)$ is set to 0.
\end{itemize}
When computing $h$ and $v$ from characteristic variables, if $0 \geq Q_+ + Q_- \propto c$, then we also set $h=0$, $v=0$. The initialization of point values is also regularized, not allowing to initialize negative water heights and replacing them by zero. Generally, wherever $h=0$, the momentum $m$ is set to zero as well.

\end{document}